\documentclass{gtart}

\def\ifplaintex{\expandafter\ifx\csname documentclass\endcsname\relax}

\def\gtp{{\mathsurround=0pt\it $\cal G\mskip-2mu$eometry \&\ 
$\cal T\!\!$opology $\cal P\!$ublications}}  

\def\recd{{\small Received:\qua\receiveddate\ifx\reviseddate\relax
\else\qquad Revised:\qua\reviseddate\fi\par}} 

\def\lognumber#1{\def\thelognumber{#1}}
\def\volumenumber#1{\def\thevolumenumber{#1}}
\def\volumeyear#1{\def\thevolumeyear{#1}}
\def\papernumber#1{\def\thepapernumber{#1}}
\def\pagenumbers#1#2{\def\startpage{#1}\def\finishpage{#2}}
\def\published#1{\def\publishdate{#1}}

\def\received#1{\def\receiveddate{#1}}
\def\revised#1{\def\reviseddate{#1}}
\def\accepted#1{\def\accepteddate{#1}}
\def\asciititle#1{\def\theasciititle{#1}}
\def\covertitle#1{\def\thecovertitle{#1}}

\let\\\par\let\thelognumber\relax\let\thevolumenumber\relax
\let\thepapernumber\relax\let\thevolumeyear\relax\let\startpage\relax
\let\finishpage\relax\let\publishdate\relax\let\receiveddate\relax
\let\reviseddate\relax\let\accepteddate\relax\let\theasciititle\relax
\let\thecovertitle\relax\let\theasciiauthors\relax
\let\theasciiabstract\relax

\let\theasciiemail\relax

\ifplaintex
\font\logobig=cmssbx10 scaled 3836
\font\logomed=cmssbx10 scaled 2557
\else
\font\logobig=cmssbx10 scaled 4200
\font\logomed=cmssbx10 scaled 2800
\fi

\long\def\makeagttitle{   
\count0=\startpage
\agt\hfill      
\hbox to 45truept{\vbox to 0pt{\vglue -13truept{\logomed A\kern -.37em{\logobig 
T}\kern -.38em G}\vss}\hss}
\break
{\small Volume \thevolumenumber\ (\thevolumeyear)
\startpage--\finishpage\nl
Published: \publishdate}

\vglue .25truein

{\parskip=0pt\leftskip 0pt plus
1fil\def\\{\par\smallskip}{\Large\bf\thetitle}\par\medskip} \vglue
0.05truein

{\parskip=0pt\leftskip 0pt plus 1fil\def\\{\par}{\sc\theauthors}
\par\medskip}%
 
\vglue 0.03truein 

{\small\leftskip 25truept\rightskip 25truept{\bf Abstract}\stdspace\theabstract

{\bf AMS Classification}\stdspace\theprimaryclass
\ifx\thesecondaryclass\relax\else; \thesecondaryclass\fi\par
{\bf Keywords}\stdspace \thekeywords\par}\vglue 7truept

}   

\ifplaintex
\font\phead=cmsl9 scaled 950
\font\pnum=cmbx10 scaled 913
\font\pfoot=cmsl9 scaled 950
\headline{\vbox to 0pt{\vskip -4.5mm\line{\small\phead\ifnum
\count0=\startpage ISSN 1472-2739 (on-line) 1472-2747 (printed)
\hfill {\pnum\folio}\else\ifodd\count0\def\\{ }%
\ifx\theshorttitle\relax\thetitle\else\theshorttitle\fi\hfill{\pnum\folio}
\else\def\\{ and }{\pnum\folio}\hfill\ifx\theshortauthors\relax\theauthors
\else\theshortauthors\fi\fi\fi}\vss}}
\footline{\vbox to 0pt{\vglue 0mm\line{\small\pfoot\ifnum\count0=\startpage
\copyright\ \gtp\hfill\else
\agt, Volume \thevolumenumber\ (\thevolumeyear)\hfill\fi}\vss}}
\else
\font=cmsl9 scaled 1050
\font\lnum=cmbx10 
\font=cmsl9 scaled 1050
\makeatletter
\def\@oddhead{{\small\lhead\ifnum\count0=\startpage ISSN 1472-2739 
(on-line) 1472-2747 (printed)\hfill {\lnum\number\count0}\else\ifodd\count0
\def\\{ }\ifx\theshorttitle\relax \thetitle \else\theshorttitle\fi\hfill
{\lnum\number\count0}\else\def\\{ and }{\lnum\number\count0}
\hfill\ifx\theshortauthors\relax 
\theauthors\else\theshortauthors\fi\fi\fi}}\def\@evenhead{\@oddhead}
\def\@oddfoot{\small\lfoot\ifnum\count0=\startpage\copyright\ \gtp\hfill\else
\agt, Volume \thevolumenumber\ (\thevolumeyear)\hfill\fi}
\def\@evenfoot{\@oddfoot}
\makeatother
\fi
\let\maketitlepage\makeagttitle
\let\makeshorttitle\maketitlepage
\let\maketitle\maketitlepage

\newwrite\gtoutfile
\long\gdef\makeheadfile{  
{\def\\{, }\def\s{ }
\immediate\openout\gtoutfile head.xxx
\immediate\write\gtoutfile{To: math@arxiv.org}
\immediate\write\gtoutfile{Subject: put OR rep NNNNN:ppppp}
\immediate\write\gtoutfile{--text follows this line--}
\immediate\write\gtoutfile{Proxy-for: \ifx\theasciiauthors\relax
\theauthors\else\theasciiauthors\fi\s<\ifx\theasciiemail\relax\theemail\else\theasciiemail\fi>}
\immediate\write\gtoutfile{\noexpand\\}
\immediate\write\gtoutfile{Authors: \ifx\theasciiauthors\relax
\theauthors\else\theasciiauthors\fi}
{\def\\{ }\immediate\write\gtoutfile{Title: \ifx\theasciititle\relax
\thetitle\else\theasciititle\fi}}
\immediate\write\gtoutfile{Subj-class: GT or SG, GR etc}
\immediate\write\gtoutfile{MSC-class: \theprimaryclass\ifx\thesecondaryclass\relax\else, \thesecondaryclass\fi}
\immediate\write\gtoutfile{Journal-ref: Algebr. Geom. Topol. \thevolumenumber\s
(\thevolumeyear) \startpage-\finishpage}
\immediate\write\gtoutfile{Comments: Published by Algebraic and
Geometric Topology at}
\immediate\write\gtoutfile{\s\s\s  http://www.maths.warwick.ac.uk/agt/AGTVol\thevolumenumber/agt-\thevolumenumber-\thepapernumber.abs.html}
\immediate\write\gtoutfile{\noexpand\\}
\immediate\write\gtoutfile{}
\ifx\theasciiabstract\relax
\immediate\write\gtoutfile{\theabstract}\else
\immediate\write\gtoutfile{\theasciiabstract}\fi
\immediate\write\gtoutfile{}
\immediate\write\gtoutfile{\noexpand\\}
\immediate\write\gtoutfile{}
\immediate\closeout\gtoutfile}}  

\def\maketitlepage{\makeagttitle\makeheadfile}
\let\makeshorttitle\maketitlepage
\let\maketitle\maketitlepage

\def\ifplaintex{\expandafter\ifx\csname documentclass\endcsname\relax}

\def\gtp{{\mathsurround=0pt\it $\cal G\mskip-2mu$eometry \&\ 
$\cal T\!\!$opology $\cal P\!$ublications}}  

\def\recd{{\small Received:\qua\receiveddate\ifx\reviseddate\relax
\else\qquad Revised:\qua\reviseddate\fi\par}} 

\def\lognumber#1{\def\thelognumber{#1}}
\def\volumenumber#1{\def\thevolumenumber{#1}}
\def\volumeyear#1{\def\thevolumeyear{#1}}
\def\papernumber#1{\def\thepapernumber{#1}}
\def\pagenumbers#1#2{\def\startpage{#1}\def\finishpage{#2}}
\def\published#1{\def\publishdate{#1}}

\def\received#1{\def\receiveddate{#1}}
\def\revised#1{\def\reviseddate{#1}}
\def\accepted#1{\def\accepteddate{#1}}
\def\asciititle#1{\def\theasciititle{#1}}
\def\covertitle#1{\def\thecovertitle{#1}}

\let\\\par\let\thelognumber\relax\let\thevolumenumber\relax
\let\thepapernumber\relax\let\thevolumeyear\relax\let\startpage\relax
\let\finishpage\relax\let\publishdate\relax\let\receiveddate\relax
\let\reviseddate\relax\let\accepteddate\relax\let\theasciititle\relax
\let\thecovertitle\relax\let\theasciiauthors\relax
\let\theasciiabstract\relax

\let\theasciiemail\relax

\ifplaintex
\font\logobig=cmssbx10 scaled 3836
\font\logomed=cmssbx10 scaled 2557
\else
\font\logobig=cmssbx10 scaled 4200
\font\logomed=cmssbx10 scaled 2800
\fi

\long\def\makeagttitle{   
\count0=\startpage
\agt\hfill      
\hbox to 45truept{\vbox to 0pt{\vglue -13truept{\logomed A\kern -.37em{\logobig 
T}\kern -.38em G}\vss}\hss}
\break
{\small Volume \thevolumenumber\ (\thevolumeyear)
\startpage--\finishpage\nl
Published: \publishdate}

\vglue .25truein

{\parskip=0pt\leftskip 0pt plus
1fil\def\\{\par\smallskip}{\Large\bf\thetitle}\par\medskip} \vglue
0.05truein

{\parskip=0pt\leftskip 0pt plus 1fil\def\\{\par}{\sc\theauthors}
\par\medskip}%
 
\vglue 0.03truein 

{\small\leftskip 25truept\rightskip 25truept{\bf Abstract}\stdspace\theabstract

{\bf AMS Classification}\stdspace\theprimaryclass
\ifx\thesecondaryclass\relax\else; \thesecondaryclass\fi\par
{\bf Keywords}\stdspace \thekeywords\par}\vglue 7truept

}   

\ifplaintex
\font\phead=cmsl9 scaled 950
\font\pnum=cmbx10 scaled 913
\font\pfoot=cmsl9 scaled 950
\headline{\vbox to 0pt{\vskip -4.5mm\line{\small\phead\ifnum
\count0=\startpage ISSN 1472-2739 (on-line) 1472-2747 (printed)
\hfill {\pnum\folio}\else\ifodd\count0\def\\{ }%
\ifx\theshorttitle\relax\thetitle\else\theshorttitle\fi\hfill{\pnum\folio}
\else\def\\{ and }{\pnum\folio}\hfill\ifx\theshortauthors\relax\theauthors
\else\theshortauthors\fi\fi\fi}\vss}}
\footline{\vbox to 0pt{\vglue 0mm\line{\small\pfoot\ifnum\count0=\startpage
\copyright\ \gtp\hfill\else
\agt, Volume \thevolumenumber\ (\thevolumeyear)\hfill\fi}\vss}}
\else
\font=cmsl9 scaled 1050
\font\lnum=cmbx10 
\font=cmsl9 scaled 1050
\makeatletter
\def\@oddhead{{\small\lhead\ifnum\count0=\startpage ISSN 1472-2739 
(on-line) 1472-2747 (printed)\hfill {\lnum\number\count0}\else\ifodd\count0
\def\\{ }\ifx\theshorttitle\relax \thetitle \else\theshorttitle\fi\hfill
{\lnum\number\count0}\else\def\\{ and }{\lnum\number\count0}
\hfill\ifx\theshortauthors\relax 
\theauthors\else\theshortauthors\fi\fi\fi}}\def\@evenhead{\@oddhead}
\def\@oddfoot{\small\lfoot\ifnum\count0=\startpage\copyright\ \gtp\hfill\else
\agt, Volume \thevolumenumber\ (\thevolumeyear)\hfill\fi}
\def\@evenfoot{\@oddfoot}
\makeatother
\fi
\let\maketitlepage\makeagttitle
\let\makeshorttitle\maketitlepage
\let\maketitle\maketitlepage

\newwrite\gtoutfile
\long\gdef\makeheadfile{  
{\def\\{, }\def\s{ }
\immediate\openout\gtoutfile head.xxx
\immediate\write\gtoutfile{To: math@arxiv.org}
\immediate\write\gtoutfile{Subject: put OR rep NNNNN:ppppp}
\immediate\write\gtoutfile{--text follows this line--}
\immediate\write\gtoutfile{Proxy-for: \ifx\theasciiauthors\relax
\theauthors\else\theasciiauthors\fi\s<\ifx\theasciiemail\relax\theemail\else\theasciiemail\fi>}
\immediate\write\gtoutfile{\noexpand\\}
\immediate\write\gtoutfile{Authors: \ifx\theasciiauthors\relax
\theauthors\else\theasciiauthors\fi}
{\def\\{ }\immediate\write\gtoutfile{Title: \ifx\theasciititle\relax
\thetitle\else\theasciititle\fi}}
\immediate\write\gtoutfile{Subj-class: GT or SG, GR etc}
\immediate\write\gtoutfile{MSC-class: \theprimaryclass\ifx\thesecondaryclass\relax\else, \thesecondaryclass\fi}
\immediate\write\gtoutfile{Journal-ref: Algebr. Geom. Topol. \thevolumenumber\s
(\thevolumeyear) \startpage-\finishpage}
\immediate\write\gtoutfile{Comments: Published by Algebraic and
Geometric Topology at}
\immediate\write\gtoutfile{\s\s\s  http://www.maths.warwick.ac.uk/agt/AGTVol\thevolumenumber/agt-\thevolumenumber-\thepapernumber.abs.html}
\immediate\write\gtoutfile{\noexpand\\}
\immediate\write\gtoutfile{}
\ifx\theasciiabstract\relax
\immediate\write\gtoutfile{\theabstract}\else
\immediate\write\gtoutfile{\theasciiabstract}\fi
\immediate\write\gtoutfile{}
\immediate\write\gtoutfile{\noexpand\\}
\immediate\write\gtoutfile{}
\immediate\closeout\gtoutfile}}  

\def\maketitlepage{\makeagttitle\makeheadfile}
\let\makeshorttitle\maketitlepage
\let\maketitle\maketitlepage

\def\ifplaintex{\expandafter\ifx\csname documentclass\endcsname\relax}

\def\gtp{{\mathsurround=0pt\it $\cal G\mskip-2mu$eometry \&\ 
$\cal T\!\!$opology $\cal P\!$ublications}}  

\def\recd{{\small Received:\qua\receiveddate\ifx\reviseddate\relax
\else\qquad Revised:\qua\reviseddate\fi\par}} 

\def\lognumber#1{\def\thelognumber{#1}}
\def\volumenumber#1{\def\thevolumenumber{#1}}
\def\volumeyear#1{\def\thevolumeyear{#1}}
\def\papernumber#1{\def\thepapernumber{#1}}
\def\pagenumbers#1#2{\def\startpage{#1}\def\finishpage{#2}}
\def\published#1{\def\publishdate{#1}}

\def\received#1{\def\receiveddate{#1}}
\def\revised#1{\def\reviseddate{#1}}
\def\accepted#1{\def\accepteddate{#1}}
\def\asciititle#1{\def\theasciititle{#1}}
\def\covertitle#1{\def\thecovertitle{#1}}

\let\\\par\let\thelognumber\relax\let\thevolumenumber\relax
\let\thepapernumber\relax\let\thevolumeyear\relax\let\startpage\relax
\let\finishpage\relax\let\publishdate\relax\let\receiveddate\relax
\let\reviseddate\relax\let\accepteddate\relax\let\theasciititle\relax
\let\thecovertitle\relax\let\theasciiauthors\relax
\let\theasciiabstract\relax

\let\theasciiemail\relax

\ifplaintex
\font\logobig=cmssbx10 scaled 3836
\font\logomed=cmssbx10 scaled 2557
\else
\font\logobig=cmssbx10 scaled 4200
\font\logomed=cmssbx10 scaled 2800
\fi

\long\def\makeagttitle{   
\count0=\startpage
\agt\hfill      
\hbox to 45truept{\vbox to 0pt{\vglue -13truept{\logomed A\kern -.37em{\logobig 
T}\kern -.38em G}\vss}\hss}
\break
{\small Volume \thevolumenumber\ (\thevolumeyear)
\startpage--\finishpage\nl
Published: \publishdate}

\vglue .25truein

{\parskip=0pt\leftskip 0pt plus
1fil\def\\{\par\smallskip}{\Large\bf\thetitle}\par\medskip} \vglue
0.05truein

{\parskip=0pt\leftskip 0pt plus 1fil\def\\{\par}{\sc\theauthors}
\par\medskip}%
 
\vglue 0.03truein 

{\small\leftskip 25truept\rightskip 25truept{\bf Abstract}\stdspace\theabstract

{\bf AMS Classification}\stdspace\theprimaryclass
\ifx\thesecondaryclass\relax\else; \thesecondaryclass\fi\par
{\bf Keywords}\stdspace \thekeywords\par}\vglue 7truept

}   

\ifplaintex
\font\phead=cmsl9 scaled 950
\font\pnum=cmbx10 scaled 913
\font\pfoot=cmsl9 scaled 950
\headline{\vbox to 0pt{\vskip -4.5mm\line{\small\phead\ifnum
\count0=\startpage ISSN 1472-2739 (on-line) 1472-2747 (printed)
\hfill {\pnum\folio}\else\ifodd\count0\def\\{ }%
\ifx\theshorttitle\relax\thetitle\else\theshorttitle\fi\hfill{\pnum\folio}
\else\def\\{ and }{\pnum\folio}\hfill\ifx\theshortauthors\relax\theauthors
\else\theshortauthors\fi\fi\fi}\vss}}
\footline{\vbox to 0pt{\vglue 0mm\line{\small\pfoot\ifnum\count0=\startpage
\copyright\ \gtp\hfill\else
\agt, Volume \thevolumenumber\ (\thevolumeyear)\hfill\fi}\vss}}
\else
\font=cmsl9 scaled 1050
\font\lnum=cmbx10 
\font=cmsl9 scaled 1050
\makeatletter
\def\@oddhead{{\small\lhead\ifnum\count0=\startpage ISSN 1472-2739 
(on-line) 1472-2747 (printed)\hfill {\lnum\number\count0}\else\ifodd\count0
\def\\{ }\ifx\theshorttitle\relax \thetitle \else\theshorttitle\fi\hfill
{\lnum\number\count0}\else\def\\{ and }{\lnum\number\count0}
\hfill\ifx\theshortauthors\relax 
\theauthors\else\theshortauthors\fi\fi\fi}}\def\@evenhead{\@oddhead}
\def\@oddfoot{\small\lfoot\ifnum\count0=\startpage\copyright\ \gtp\hfill\else
\agt, Volume \thevolumenumber\ (\thevolumeyear)\hfill\fi}
\def\@evenfoot{\@oddfoot}
\makeatother
\fi
\let\maketitlepage\makeagttitle
\let\makeshorttitle\maketitlepage
\let\maketitle\maketitlepage

\newwrite\gtoutfile
\long\gdef\makeheadfile{  
{\def\\{, }\def\s{ }
\immediate\openout\gtoutfile head.xxx
\immediate\write\gtoutfile{To: math@arxiv.org}
\immediate\write\gtoutfile{Subject: put OR rep NNNNN:ppppp}
\immediate\write\gtoutfile{--text follows this line--}
\immediate\write\gtoutfile{Proxy-for: \ifx\theasciiauthors\relax
\theauthors\else\theasciiauthors\fi\s<\ifx\theasciiemail\relax\theemail\else\theasciiemail\fi>}
\immediate\write\gtoutfile{\noexpand\\}
\immediate\write\gtoutfile{Authors: \ifx\theasciiauthors\relax
\theauthors\else\theasciiauthors\fi}
{\def\\{ }\immediate\write\gtoutfile{Title: \ifx\theasciititle\relax
\thetitle\else\theasciititle\fi}}
\immediate\write\gtoutfile{Subj-class: GT or SG, GR etc}
\immediate\write\gtoutfile{MSC-class: \theprimaryclass\ifx\thesecondaryclass\relax\else, \thesecondaryclass\fi}
\immediate\write\gtoutfile{Journal-ref: Algebr. Geom. Topol. \thevolumenumber\s
(\thevolumeyear) \startpage-\finishpage}
\immediate\write\gtoutfile{Comments: Published by Algebraic and
Geometric Topology at}
\immediate\write\gtoutfile{\s\s\s  http://www.maths.warwick.ac.uk/agt/AGTVol\thevolumenumber/agt-\thevolumenumber-\thepapernumber.abs.html}
\immediate\write\gtoutfile{\noexpand\\}
\immediate\write\gtoutfile{}
\ifx\theasciiabstract\relax
\immediate\write\gtoutfile{\theabstract}\else
\immediate\write\gtoutfile{\theasciiabstract}\fi
\immediate\write\gtoutfile{}
\immediate\write\gtoutfile{\noexpand\\}
\immediate\write\gtoutfile{}
\immediate\closeout\gtoutfile}}  

\def\maketitlepage{\makeagttitle\makeheadfile}
\let\makeshorttitle\maketitlepage
\let\maketitle\maketitlepage

\def\ifplaintex{\expandafter\ifx\csname documentclass\endcsname\relax}

\def\gtp{{\mathsurround=0pt\it $\cal G\mskip-2mu$eometry \&\ 
$\cal T\!\!$opology $\cal P\!$ublications}}  

\def\recd{{\small Received:\qua\receiveddate\ifx\reviseddate\relax
\else\qquad Revised:\qua\reviseddate\fi\par}} 

\def\lognumber#1{\def\thelognumber{#1}}
\def\volumenumber#1{\def\thevolumenumber{#1}}
\def\volumeyear#1{\def\thevolumeyear{#1}}
\def\papernumber#1{\def\thepapernumber{#1}}
\def\pagenumbers#1#2{\def\startpage{#1}\def\finishpage{#2}}
\def\published#1{\def\publishdate{#1}}

\def\received#1{\def\receiveddate{#1}}
\def\revised#1{\def\reviseddate{#1}}
\def\accepted#1{\def\accepteddate{#1}}
\def\asciititle#1{\def\theasciititle{#1}}
\def\covertitle#1{\def\thecovertitle{#1}}

\let\\\par\let\thelognumber\relax\let\thevolumenumber\relax
\let\thepapernumber\relax\let\thevolumeyear\relax\let\startpage\relax
\let\finishpage\relax\let\publishdate\relax\let\receiveddate\relax
\let\reviseddate\relax\let\accepteddate\relax\let\theasciititle\relax
\let\thecovertitle\relax\let\theasciiauthors\relax
\let\theasciiabstract\relax

\let\theasciiemail\relax

\ifplaintex
\font\logobig=cmssbx10 scaled 3836
\font\logomed=cmssbx10 scaled 2557
\else
\font\logobig=cmssbx10 scaled 4200
\font\logomed=cmssbx10 scaled 2800
\fi

\long\def\makeagttitle{   
\count0=\startpage
\agt\hfill      
\hbox to 45truept{\vbox to 0pt{\vglue -13truept{\logomed A\kern -.37em{\logobig 
T}\kern -.38em G}\vss}\hss}
\break
{\small Volume \thevolumenumber\ (\thevolumeyear)
\startpage--\finishpage\nl
Published: \publishdate}

\vglue .25truein

{\parskip=0pt\leftskip 0pt plus
1fil\def\\{\par\smallskip}{\Large\bf\thetitle}\par\medskip} \vglue
0.05truein

{\parskip=0pt\leftskip 0pt plus 1fil\def\\{\par}{\sc\theauthors}
\par\medskip}%
 
\vglue 0.03truein 

{\small\leftskip 25truept\rightskip 25truept{\bf Abstract}\stdspace\theabstract

{\bf AMS Classification}\stdspace\theprimaryclass
\ifx\thesecondaryclass\relax\else; \thesecondaryclass\fi\par
{\bf Keywords}\stdspace \thekeywords\par}\vglue 7truept

}   

\ifplaintex
\font\phead=cmsl9 scaled 950
\font\pnum=cmbx10 scaled 913
\font\pfoot=cmsl9 scaled 950
\headline{\vbox to 0pt{\vskip -4.5mm\line{\small\phead\ifnum
\count0=\startpage ISSN 1472-2739 (on-line) 1472-2747 (printed)
\hfill {\pnum\folio}\else\ifodd\count0\def\\{ }%
\ifx\theshorttitle\relax\thetitle\else\theshorttitle\fi\hfill{\pnum\folio}
\else\def\\{ and }{\pnum\folio}\hfill\ifx\theshortauthors\relax\theauthors
\else\theshortauthors\fi\fi\fi}\vss}}
\footline{\vbox to 0pt{\vglue 0mm\line{\small\pfoot\ifnum\count0=\startpage
\copyright\ \gtp\hfill\else
\agt, Volume \thevolumenumber\ (\thevolumeyear)\hfill\fi}\vss}}
\else
\font=cmsl9 scaled 1050
\font\lnum=cmbx10 
\font=cmsl9 scaled 1050
\makeatletter
\def\@oddhead{{\small\lhead\ifnum\count0=\startpage ISSN 1472-2739 
(on-line) 1472-2747 (printed)\hfill {\lnum\number\count0}\else\ifodd\count0
\def\\{ }\ifx\theshorttitle\relax \thetitle \else\theshorttitle\fi\hfill
{\lnum\number\count0}\else\def\\{ and }{\lnum\number\count0}
\hfill\ifx\theshortauthors\relax 
\theauthors\else\theshortauthors\fi\fi\fi}}\def\@evenhead{\@oddhead}
\def\@oddfoot{\small\lfoot\ifnum\count0=\startpage\copyright\ \gtp\hfill\else
\agt, Volume \thevolumenumber\ (\thevolumeyear)\hfill\fi}
\def\@evenfoot{\@oddfoot}
\makeatother
\fi
\let\maketitlepage\makeagttitle
\let\makeshorttitle\maketitlepage
\let\maketitle\maketitlepage

\newwrite\gtoutfile
\long\gdef\makeheadfile{  
{\def\\{, }\def\s{ }
\immediate\openout\gtoutfile head.xxx
\immediate\write\gtoutfile{To: math@arxiv.org}
\immediate\write\gtoutfile{Subject: put OR rep NNNNN:ppppp}
\immediate\write\gtoutfile{--text follows this line--}
\immediate\write\gtoutfile{Proxy-for: \ifx\theasciiauthors\relax
\theauthors\else\theasciiauthors\fi\s<\ifx\theasciiemail\relax\theemail\else\theasciiemail\fi>}
\immediate\write\gtoutfile{\noexpand\\}
\immediate\write\gtoutfile{Authors: \ifx\theasciiauthors\relax
\theauthors\else\theasciiauthors\fi}
{\def\\{ }\immediate\write\gtoutfile{Title: \ifx\theasciititle\relax
\thetitle\else\theasciititle\fi}}
\immediate\write\gtoutfile{Subj-class: GT or SG, GR etc}
\immediate\write\gtoutfile{MSC-class: \theprimaryclass\ifx\thesecondaryclass\relax\else, \thesecondaryclass\fi}
\immediate\write\gtoutfile{Journal-ref: Algebr. Geom. Topol. \thevolumenumber\s
(\thevolumeyear) \startpage-\finishpage}
\immediate\write\gtoutfile{Comments: Published by Algebraic and
Geometric Topology at}
\immediate\write\gtoutfile{\s\s\s  http://www.maths.warwick.ac.uk/agt/AGTVol\thevolumenumber/agt-\thevolumenumber-\thepapernumber.abs.html}
\immediate\write\gtoutfile{\noexpand\\}
\immediate\write\gtoutfile{}
\ifx\theasciiabstract\relax
\immediate\write\gtoutfile{\theabstract}\else
\immediate\write\gtoutfile{\theasciiabstract}\fi
\immediate\write\gtoutfile{}
\immediate\write\gtoutfile{\noexpand\\}
\immediate\write\gtoutfile{}
\immediate\closeout\gtoutfile}}  

\def\maketitlepage{\makeagttitle\makeheadfile}
\let\makeshorttitle\maketitlepage
\let\maketitle\maketitlepage

\lognumber{40}
\volumenumber{2}
\volumeyear{2002}
\papernumber{40}
\pagenumbers{1001}{1050}
\received{30 November 2002}
\revised{15 July 2002}
\accepted{21 July 2002}
\published{5 November 2002}

\usepackage{epsf, amssymb}
\def\codim{\mathop{\rm codim}}

\theoremstyle{plain}
\newtheorem{prop}[equation]{Proposition}

\newtheorem{lemma}[equation]{Lemma}
\theoremstyle{definition}
\newtheorem{remark}[equation]{Remark}
\newtheorem{defin}[equation]{Definition}
\newtheorem{nota}[equation]{Notation}
\newtheorem{note}[equation]{Notation}
\catcode`\@=11
\def\eqalign#1{\null \,\vcenter {\openup \jot \m@th \ialign 
{\strut \hfil $\displaystyle {##}$&$\displaystyle {{}##}$\hfil \crcr 
#1\crcr }}\,}\catcode`\@=12
\font\tensym=msbm10   
\font\sevensym=msbm7
\font\fivesym=msbm5
\newfam\symfam
\textfont\symfam=\tensym
\scriptfont\symfam=\sevensym
\scriptscriptfont\symfam=\fivesym
\def\sym{\fam\symfam\tensym}

\def\R{{\sym R}}

\def\N{{\sym N}}
\def\Z{{\sym Z}}

\def\qs{\forall\,}

\def\imp{\ \hbox{$\Relbar\mkern -10mu\Rightarrow$~}}

\def\dm{{\cal D}(M)}
\def\dn{{\cal D}_n(M)}
\def\b{$\bullet$\ }
\def\f{f}
\let\0=\emptyset
\def\diagram#1{\def\normalbaselines {\lineskip=5pt\baselineskip=0pt
\lineskiplimit=1pt}\matrix{#1}}
\def\hfl#1{\smash{\mathop{\hbox to 10mm{\rightarrowfill}}\limits^{\textstyle
#1}}}

\def\vf{\left\downarrow\vbox to 5mm{}\right.}
\def\A{{\cal A}}

\def\g{\Gamma}
\def\ed{{\cal D}}
\def\edm{\ed_n(M)}
\def\ew{{\cal D}^W}

\def\dl{D_{n,k}(M)}
\def\dk{D_{n,k}}
\def\dw{D^W_{n,k}}
\def\ddw{D^W_N}
\def\eg{{\cal D}^g}
\let\gu=\gamma
\def\egn{{\eg_{n,k}}}
\def\dg{D^g_{n,k}}

\def\cc{{C'_\g}}
\def\wg{{\cal W}_m(\gu)}
\def\ank{\A_n^k(M)}
\def\wgg{{\cal W}_{m(\gu)}(\gu)}
\def\ce{{\cal E}_s}
\def\ee{{\cal E}}

\def\wx{W_s}
\let \inc=\subseteq
\let\O=\Omega
\def\h{{\cal H}(G)}
\let\ph=\varphi
\def\cb{\overline \c}
\def\gr{\overline C(G)}
\def\poids{\beta}
\def\ano{\alpha}

\def\ss{\co(G,S)}

\def\sd{S^2}
\let \ps=\psi
\def\rt{\R^3}
\def\dd{D^2}
\def\D{l_s}
\def\lo{l_z}
\def\al{Q}

\def\s{\sigma}
\let\om=\omega
\let\th=\theta
\def\bb{{\bar B}}
\def\enk{{\cal G}
}\def\wt{\widetilde W}
\def\ba{{\bar A}}
\def\vep#1{\vcenter{\hbox{\epsfbox{#1.eps}}}}
\def\epfb#1{\left(\,\vcenter{\hbox{\epsfbox{#1.eps}}}\,\right)}

\let\del=\partial

\def\bu{{\bar U}}
\def\ha{{\widehat A}}
\def\ra{R_A}

\def\hh{{\cal H}(\g)}
\def\rb{R_\ha}
\def\1{+\infty}
\let\lra=\longrightarrow
\def\co{C}
\def\dem{E_{1\over 2}}

\begin{document}
\title{The Configuration space integral for links in $\mathbb R^3$}
\asciititle{The Configuration space integral for links in R^3}
\covertitle{The Configuration space integral for links in $\noexpand\bf R^3$}
\author{Sylvain Poirier}
\address{3 bis rue Auguste Dollfus, 76600 Le Havre, France}
\email{spoirier@lautre.net}

\begin{abstract}
The perturbative expression of Chern-Simons theory for links in 
Euclidean 3-space is a linear combination of integrals on 
configuration spaces. This has successively been studied by 
Guadagnini, Martellini and Mintchev, Bar-Natan, Kontsevich, 
Bott and Taubes, D. Thurston, Altschuler and Freidel, 
Yang and others. We give a self-contained version
of this study with a new choice of compactification, and we formulate a 
rationality result. 
\end{abstract}

\maketitle

\section{Introduction}
\setcounter{equation}{0}
\def\c{{C_\g}}
\def\cg{C_\g^{tot}}
\def\truc{x}
The perturbative expression of Chern-Simons theory for links in 
Euclidean 3-space in the Landau gauge with the generic (universal)
gauge group, that we call 
``configuration space integral'' (according
to a suggestion of D.Bar-Natan), is a linear combination of integrals on 
configuration spaces. It has successively been studied by 
E. Guadagnini, M. Martellini and M. Mintchev \cite{gua}, 
D.Bar-Natan \cite{bn'}, M. Kontsevich \cite{ko}, 
Bott and Taubes \cite{bt}, Dylan Thurston
\cite{th} (who explained in details the notion of degree of a map),
Altschuler and Freidel \cite{af}, Yang \cite{ya} and others.
We give a self-contained version
of this study with a new choice of compactification, and formulate a 
rationality result.

\subsection{Definitions}
Let $M$ be a compact one-dimensional manifold
with boundary.
Let $L$ denote an imbedding of $M$ into $\rt$. We say that $L$ is a 
{\em link} if we moreover have the condition that the boundary of $M$ is
empty. The most important results in this article only concern the case of 
links, but we allow the general case in our definitions in order to
prepare the next article.

And a link $L$ it is a {\em knot} when $M=S^1$. 
\ppar

A {\em diagram with support} $M$, is  the data of :

(1)\qua A graph made of a finite set $V=U\cup T$ of vertices
and a set $E$ of edges which are pairs of elements of $V$, 
such that : the
elements of $U$ are univalent (belong to precisely one edge) 
and the elements of $T$ are trivalent, and
every connected component of this graph meets $U$. 

(2)\qua An isotopy class $\s$ of injections $i$
of $U$ into the interior of $M$.
\ppar

This definition of a diagram differs from the usual one in that it
excludes the double edges $\epfb{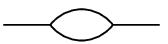}$ and the loops ($\circ\kern -1pt-$).
This will be justified by the fact that the integral that we shall define
makes no sense on loops, and by Proposition \ref{34} which will allow us to
exclude diagrams with a double edge, among others.

For a diagram $\g$, the {\em degree} of $\g$ is the
number $$\deg\g={1\over 2}\#V=\#E-\#T={1\over 3}(\#E+\#U).$$
An {\em isomorphism} between two diagrams
$\g=(V,E,\s)$ and $\g'=(V',E',\s')$, is a bijection from $V$ to $V'$
which transports $E$ to $E'$ and $\s$ to $\s'$.

Let $|\g|$ denote the number of automorphisms of $\g$.
\ppar

Let $\dm$ denote the set of diagrams with support $M$ 
up to isomorphism, and $\dn$ its subset made of 
the degree $n$ diagrams, for all $n\in\N$.

${\cal D}_0(M)$ has one element which is the empty diagram.

${\cal D}_1(S^1)$ also has only one element, which will be 
denoted by $\th$, for it looks like the letter $\th$.

An {\em orientation} $o$ of a diagram $\g$ is the data of 
an orientation of all its vertices, where the {\em orientation} of
a trivalent vertex is a cyclic order on the set of the three edges
it belongs to, and the {\em orientation} of a univalent vertex is a
local orientation of $M$ near its image by $i$.

Let $\A_n(M)$ denote the real vector space generated by the set 
of oriented diagrams of degree $n$, and quotiented
by the following AS, IHX and STU relations. 
The image of an oriented diagram $\g$ in the space $\A_n$ is denoted by $[\g]$.

The AS relation says that if $\g$ and $\g'$ differ by the 
orientation of exactly one vertex, then $[\g']=-[\g]$.

The IHX and STU relations are respectively defined by the following
formulas. 
$$\vep{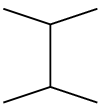}=\vep{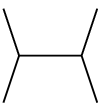}-\vep{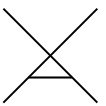}\qquad,\qquad\vep{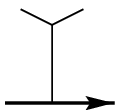}=\vep{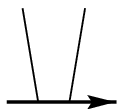}-\vep{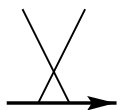}$$
These formulas relate diagrams which are identical outside
one place, where they differ according to the figures.
By convention for the figures, the orientation at each trivalent vertex is 
the given counterclockwise
order of the edges, and the orientation of each univalent vertex is given
by the drawn orientation of the bold line (which represents $M$ locally).

Now let $\A(M)$ denote the space $$\A(M)=\prod_{n\in\N}\A_n(M).$$
A {\em configuration} of a diagram $\g\in\dm$ on $L$ is a map
from the set $V$ of vertices of $\g$ to $\rt$, which is injective on 
each edge, and which coincides on $U$ with $L\circ i$ for some $i$ in 
the class $\s$. The {\em configuration space} $\c(L)$ 
(or simply $\c$ if $L$ is fixed) of a diagram $\g\in\dm$ 
on $L$ is the set
of these configurations. 

$\c$ has a natural structure of a smooth manifold of dimension 
$\#U+3\#T=2\#E$.
We shall canonically associate an orientation of $\c$ to
the data of an orientation of $\g$ and an orientation of all edges of $\g$,
such that it is reversed when we change the orientation of one vertex 
or one edge (Subsection \ref{ori}).

When orientations of all edges are chosen, we have the following canonical
map $\Psi$ from $\c$ to $(\sd)^E$: for a configuration $x\in\c$ and
$e=(v,v')$, $\Psi(x)_e$, is 
the direction of $x(v')-x(v)$.

Let $\om$ denote the standard volume form on $\sd$ with total mass 1, and let
$$\O=\bigwedge\nolimits^E\om$$denote the corresponding volume form on $(\sd)^E$.

Now the {\em configuration space integral} is $Z(L)\in \A(M)$:
$$Z(L)=\sum_{\g\in \dm} {I_L(\g)\over|\g|}[\g]$$
where
$$I_L(\g)=\int_{\c(L)}\Psi^*\O.$$
We shall see that this integral $I_L(\g)$ converges. The sign of
the product $I_L(\g)[\g]$ will not depend on choices
of orientations. Then $Z(L)$ also converges because it is a finite sum 
in each $\A_n(M)$.

By convention, the zero-degree part of $Z(L)$ is $1=[\0]$.

As an introduction to the study of $Z(L)$, let us first recall the properties
of its degree one part when $L$ is a link.

There are two types of degree one diagrams: the chords $c_{m,m'}$ between
two different components $m$ and $m'$ of $M$, and the chords $\th_m$ with
both ends on the same component $m$ of $M$. There are no IHX or STU
relations between them, so the degree one part of $Z$ consists of the
data of the integrals $I_L(c_{m,m'})$ and $I_L(\th_m)$.

The space $C_{c_{m,m'}}(L)$ is identical to $m\times m'$, so, since $M$ has
no boundary, it is diffeomorphic to the torus $S^1\times S^1$.
Now, $\Psi$ is a smooth map from $m\times m'$ to $\sd$, so the integral
$I_L(c_{m,m'})=\int_{c_{m,m'}}\Psi^*\om$ is equal to the degree of 
the map $\Psi$: it takes values in $\Z$. Precisely, it is equal 
to the linking number of $L(m)$ and $L(m')$.

The space $C_\th$ is diffeomorphic to the noncompact surface
$S^1\times ]0,1[$. So the integral
$I_L(\th)$ is not the degree of a map. In fact it is known as the {\em Gauss 
integral}, usually considered in the case of a knot. It 
takes all real values on each isotopy class of knots. In the case of an
almost planar knot, it approaches an integer which is the writhe or {\em 
self-linking number}
of this knot (the number of crossings counted with signs).

In order to express some theorems about this expression $Z(L)$, we need 
first to recall some usual algebraic structures.

\subsection{Algebraic structures on the spaces $\A$ of diagrams}

For disjoint $M$ and $M'$ included in some $M''$ such that their interiors 
do not meet, one can consider the bilinear map from 
$\A(M)\times\A(M')$ to $ \A(M'')$
defined by: for every oriented diagrams $\g$ with support $M$
and $\g'$ with support $M'$, $[\g]\cdot[\g']=[\g\sqcup\g']$. This
bilinear map is graded, that is, it sends $\A_n(M)\times\A_{n'}(M')$
to $\A_{n+n'}(M'')$.

Now, let $J=[0,1]$. It is well-known \cite{bn} that the bilinear product on $\A(J)$
defined by gluing the two copies of $J$ preserving the orientation, defines
a commutative algebra structure with unit $1=[\0]$, 
and that $\A(J)$ can be identified with $\A(S^1)$.

Similarly, for every $M$ and every connected component $m$ of $M$ with
an orientation (either with boundary or not), the insertion of $J$ at
some point of $m$ preserving the orientation provides an 
$\A(J)$-module structure on $\A(M)$: it is well-known \cite{bn} 
that the result does not depend on the place of the insertion, 
modulo the AS, IHX and STU relations.

\subsection{Results on the configuration space integral}

The convergence of the integrals $I_L(\g)$ comes from the following 
proposition which will be proved in Section \ref{comp}.

\begin{prop}\label{compac} 
There exists a compactification of $\c$ into a manifold
with corners $\cb$, on which the canonical map $\Psi$ extends smoothly.
\end{prop}
\def\anom{\ano^{(m)}}

Bott and Taubes \cite{bt}
first proved the convergence of the integrals by compactifying the
configuration space $\c$ into a manifold with corners, 
where $\Psi^*\O$ extends smoothly. (Their compactification
was a little stronger than ours).

The notion of a manifold with corners is defined in their Appendix; 
briefly, it means every point has a neighbourhood diffeomorphic to 
$]-1,1[^p\times[0,1[^{n-p}$ for some $p$ from $0$ to $n$.

They restricted their interest to the case of knots, and they studied
the variations of these integrals during isotopies, in order to check the
invariance. So, they applied the Stokes theorem to express the variations
of $I_L(\g)$, in terms of integrals on the fibered spaces over the time, made
of the faces of the $\cb$. 
They found equalities between the integrals
at certain faces of different $\cb$, where the diagrams $\g$ are the terms of
an IHX or STU relation, and proved cancellation on the other faces except
those of a special kind called ``anomalous'', which follow in some way 
the tangent map of $L$. The global effect of the contributions 
of the anomalous faces to the
variations of the integral can be expressed by a 
universal constant $\ano\in \A$, the so-called {\it anomaly}.

Using many of their arguments, we shall prove the following proposition 
for links (when $\del M=\0$):

\begin{prop}\label{varia} The variations of $Z(L)$ on each isotopy class of links are expressed
by the formula $$Z(L)=Z^0(L) \prod_m \exp\left({I_L(\th_m)}\anom\right)$$
where $m$ runs over the connected components of $M$, and :

$Z^0(L)\in \A(M)$ is an isotopy invariant of $L$,\nl
$\ano$ is a universal constant in $\A(J)$ called the {\em anomaly},\nl
$\anom$ means $\ano$ acting on $\A(M)$ on the component $m$.
\end{prop}

(But if the boundary of $M$ was not empty, it would produce another boundary
of the $\cb$ and thus uncontrollable variations of $Z$.)

Altschuler and Freidel \cite{af} have written this formula
in the case of knots. They and Dylan Thurston \cite{th} also prove that 
\begin{prop}$Z^0$ is a universal 
Vassiliev invariant.\end{prop} 

The degree one part of $\ano$ is ${[\th]\over 2}$, and it is not difficult to see that
all even degree parts of $\ano$ cancel for symmetry reasons (Lemma \ref{sym2}).

We shall also prove in Section \ref{v3} that

\begin{prop} 
The degree 3 part of the anomaly vanishes.
\end{prop}

But first we shall prove the following Proposition \ref{34} 
in Section \ref{pr34}:

\begin{nota}
Let $A$ be any nonempty subset of $V=U\cup T$. We denote
$$\eqalign{E_A&=\{e\in E\,|\,e\inc A\}\cr E'_A&=\{e\in E\,|\,\#(e\cap A)=1\}.}$$
The cardinalities of these sets are related by the following formula 
which expresses the two ways of counting the half-edges of the vertices in $A$:
\begin{equation}
\label{count}
2\#E_A+\#E'_A=3\#(A\cap T)+\#(A\cap U)
\end{equation}
\end{nota}

\begin{defin}\label{3con}
A diagram $\g$ is said to be {\it principal} if for any 
$A\inc T$ such that $ \#A>1$ we have $\#E'_A\geq 4$.

A diagram $\g$ is said to be {\it subprincipal}\/ if it obeys the two following
conditions:

(1)\qua For any 
$A\inc T$ such that $\#A>1$, we have $\#E'_A\geq 3$

(2)\qua For any two $A,A'\inc T$ such that $\#A>1$, $\#A'>1$ and 
$\#E'_A=\#E'_{A'}=3$, we have $A\cap A'\not=\0$.
\end{defin}

The interest of these definitions comes from the following properties:

\begin{prop}\label{34} For a given $\g$, there exists a $L$ such
that $\codim\Psi(\c)=0$ if and only if $\g$ is principal.
And the property that $\g$ is subprincipal is a necessary condition so that 
$\codim\Psi(\c)\leq 1$. 
\end{prop}

\begin{remark} The equivalence is also probably true for the 
second case, when $L$ is generic.
But we shall not use this result in this article. In fact, as what we only
need is a necessary condition so that 
$\codim\Psi(\c)\leq 1$, we shall not make use of Condition (2) in the definition
of a subprincipal diagram.

In fact, Condition (1) implies that any sets $A,A'\inc T$ such that $\#A>1$, 
$\#A'>1$ and $\#E'_A=\#E'_{A'}=3$ are either disjoint or one is included in 
the other (this is a bit tedious to prove and we shall not need it).
\end{remark}

The idea of the proof of Proposition \ref{varia} is that for every $n$, 
the degree $n$ part $Z^0_n$ of $Z^0$ can be seen as the degree of a certain map 
from a glued union of configuration spaces which are given rational multiples
of the classes of their diagrams in $\A(M)$ as coefficients, 
to a certain manifold $(\sd)^N$ (where we can take for $N$ the maximum
number of edges of degree $n$ diagrams).

This implies that $Z^0_n$ is a rational linear combination of the $[\g]$
for the degree $n$ diagrams $\g$. More precisely we have the following proposition, which will be proved together with Proposition \ref{varia} in Section \ref{secano}.

\begin{nota}
For a diagram $\g$, let $u_\g$ be the number of univalent vertices of $\g$, and $e_\g$ be
its number of edges (we have $u_\g+e_\g=3\deg\g)$.

Let $k\in\Z$, $k\leq 2n$, and let $\ank$ be the quotient space of 
$\A_n(M)$ by the subspace generated by the subprincipal
diagrams $\g$ such that $u_\g=k-1$. (This implies, through the STU relation,
> the cancellation of all subprincipal diagrams with $u_\g<k$.)

We must suppose that $k\leq 2n$, else we would have $\ank=\{0\}$.
Let $Z_n^k$ be the image of $Z_n$ in $\ank$. 
Proposition \ref{34} implies that $\A_n^2(M)=\A_n(M)$, and that 
$\A_n^3(M)=\A_n(M)$ if $n>1$.
\end{nota}

The reason for this construction of $\ank$ will appear in the proof of
Lemma \ref{rstu} which is used to obtain the following rationality result:

\begin{prop}
\label{mult}
We suppose $L$ is such that for all $m$, $I_L(\th_m)$ is an integer. 
Let $n\geq1$, $k\in\Z$ and $N=3n-k$. 
Then $Z_n^k$ belongs to the integral lattice
generated by the elements$${(N-e_\g)!\over N!\;2^{e_\g}}[\g]$$
in the space $\ank$, where $\g$ runs over the degree $n$ principal
diagrams.
\end{prop}

\begin{remark}\label{improve}
Here are two results which will be proved in a future article:

\b\qua The above proposition can be improved by saying that $Z_n$ 
belongs to the integral lattice
generated by the elements$${[\g]\over e_\g!\;2^{e_\g}}$$
where $\g$ runs over the set of principal diagrams.

\b\qua The degree 5 part of the anomaly cancels (this uses a Maple program).
\end{remark}

The present rationality result can be compared to the one obtained 
for the Kontsevich integral: it was proved in \cite{le} that the 
degree $n$ part of the Kontsevich integral is 
a linear combination of chord diagrams (or of all diagrams, 
which is the same), where the
coefficients are multiples of$${1\over (1!2!\cdots n!)^4(n+1)}.$$
For large values of $n$, these denominators are greater than
ours, but have smaller prime factors.

\ppar

I would like to express my thanks to Christine Lescop for the help
during the preparation of this paper. 
I also thank Gregor Masbaum, Pierre Vogel, Alexis Marin, Su-Win Yang, 
Dylan Thurston, Dror Bar-Natan, Simon Willerton and 
Michael Polyak for their advice and comments.

\section{Compactification}\label{comp}
\setcounter{equation}{0}

\subsection{The compactified space $\h$ of a configuration space $\co(G)$ of a graph}
\let\gr=\h

\begin{nota} If $A$ is a finite set with at least two elements, let
$C^A$ denote the space
of nonconstant maps from 
$A$ to $\rt$ quotiented by the translation-dilations group 
(that is the group of translations and positive homotheties of $\rt$).
\end{nota}

Note that the space $C^A$ is diffeomorphic to the sphere $S^{3\#A-4}$: choose
an element $x\in A$ to be at the origin; then the images of the other elements
of $A$ are to be considered modulo the dilations with center the origin.
Thus it is compact.

\begin{nota}\label{notgra} Let $G$ be a graph defined as a pair $G=(V,E)$ 
where $V$ is the set of ``vertices'' and $E$
is the set of ``edges'': $V$ is a finite set and $E$ is a set of 
pairs of elements of $V$. We suppose
that $\#V\geq 2$ and that $G$ is connected. 

We define the configuration space $\co(G)$ as the space of maps
from $V$ to $\rt$ which map the two ends of every edge to two different points, modulo
the dilations and translations (it is a dense open subset of $C^V$).

By an abuse of notation, any subset $A$ of $V$ will also mean the graph
with the set $A$ of vertices and the set $E_A$ of edges (the set of edges 
$e\in E$ such that $e\inc A$).  

Let $R$ be the set of connected subsets $A$ of $V$ such that 
$\#A\geq 2$.

Let $\h$ be the subset of $\prod_{A\in R}C^A$ made of the
$x=(x_A)_{A\in R}$ such that $\qs A,B\in R$, $ A\inc B\imp$ the restriction 
of $x_B$ to $A$ is either the constant map or the map $x_A$.
\end{nota}

We define a 
canonical imbedding of $\co(G)$ in $\gr$ by restricting any $f\in \co(G)$ 
to each of the $A\in R$: indeed, $\qs A\in R,\ \qs f\in \co(G)$, $f$ 
is not constant on $A$ because $A$ contains at least one edge, so the
restriction of $f$ to $A$ is a well-defined element of
$C^A$. Moreover, this is an imbedding because
$V\in R$. 

\def\L{{\cal L}}
We can see that $\h$ is compact, as a closed subset of the compact manifold 
$\prod_{A\in R}C^A$, for it is defined as an intersection of closed sets:

$\qs A,B\in R$ such that $ A\inc B$, let us see why $\{(x_A,x_B)\in
C^A\times C^B|$ the restriction of $x_B$ to $A$ is either 
the constant map or the map $x_A\}$ is closed, in the following way:
by identifying each of $C^A$, $C^B$ as a sphere of representatives a convenient
way (fixing one vertex in $A$ to the origin) in the respective linear space 
$\L$, $\L'$ with the canonical linear projection $\pi$ from $\L'$ onto 
$\L$, then the above set is the projection by 
$(x_A,x_B,\lambda)\mapsto (x_A,x_B)$
of the closed thus compact set of $(x_A,x_B,\lambda)\in C^A\times C^B\times
[0,1]$ such that $\pi(x_B)=\lambda x_A$.

In the next subsection, we shall prove that $\h$ is a manifold
with corners, and that it is a compactification of $\co(G)$.
Now we are going to describe the strata of $\gr$. First, let us introduce
this description intuitively:

An element $x$ of a stratum of $\gr$ will be a limit of elements of $\co(G)$.
This limit is first described by its projection to $C^V$; but some edges will
collapse if $x$ is not in the interior stratum $\co(G)$. Group the collapsing
edges into connected components $A$, zoom in on each such $A$ 
and describe the relative positions of the vertices in 
$A$: it is an element of $C^A$. If some edges in $A$ still collapse, do the
same operations until all the relative positions of points in connected 
subdiagrams are described. The set of all sets $A$ on which
you have zoomed in forms a subset $S$ of $R$.

\begin{nota}\label{nots}
In all the following, $S$ will denote any subset of $R$ such that $V\in S$ and any two elements
of $S$ are either disjoint or included one into the other. Let $S'=S-\{V\}$.

For any $A\in R$, let $G/A$ denote the graph obtained by identifying the elements of $A$ into one vertex, and deleting the edges in $E_A$.

For any $A\in S$, let $A/S$ denote the graph $A$ quotiented 
by the greatest elements $A_i$ of $S$ strictly included
in $A$ (they are disjoint and any element of $S$ strictly included in $A$
is contained in one of them).

Identify the space $\co(A/S)$ to a subset of $C^A$ in the natural way (each $x\in\co(A/S)$ is 
then constant on each of the sets $A_i$) and
let$$\ss=\prod_{A\in S} \co(A/S)\inc \prod_{A\in S} C^A.$$
$\qs A\inc V, A\not=\0$, let $\ba$ denote the smallest element of $S$ containing $A$;

$\qs A\in S'$, let $\ha$ denote the smallest element of $S$ that strictly 
contains $A$.
\end{nota}

Let us first check that the last definitions make sense. 

For any nonempty subset $A$ of $V$, the smallest element of $S$ containing $A$
is well-defined because first $A\inc V\in S$, then two elements of 
$S$ which contain
$A$ cannot be disjoint, so one of them is included in the other.

The definition of $\ha$ is justified in the same way. This gives $S$ a
tree structure.

\def\bx{\tilde x}

\begin{lemma}\label{plo} 
{\rm(1)}\qua For all S and all $x\in \ss$, there is a unique $\bx\in \h$ which 
is an extension of $x$ to $R$. 
This defines an imbedding of $\ss$ into $\gr$.

{\rm(2)}\qua This $\bx$ has the property that
for all $A,B\in R$ with $A\inc B$, the restriction of $\bx_{B}$ to
$A$ is constant if and only if $\ba\not=\bar{B}$.

{\rm(3)}\qua For all $x\in\gr$, there is a unique $S$ such that the restriction of 
$x$ to $S$ belongs
to $\ss$.
\end{lemma}

\proof

(1)\qua First let us see the construction and uniqueness of $\bx$.

Since $\bx$ must belong to $\gr$, for all $A\in R$, the restriction of
$\bx_\ba=x_\ba$ to $A$ must be either constant or equal to $\bx_A$. But it
is not constant for the following reason:

By definition of $\ba$, $A$ is not reduced to one vertex in $\ba/S$.
Furthermore, since $A$ is connected, it has at least one common edge 
with $\ba/S$.
Then, since $x_\ba\in\co(\ba/S)$, it cannot be constant on this edge.

Thus, for all $A\in R$, the value of $\bx_A$ is determined by $x_\ba$. 
This proves the uniqueness of $\bx$. 

As for the existence of $\bx$, the fact that this $\bx$ constructed above 
is actually an element of $\h$, will be a consequence of (2). 

Now, the fact that this defines an imbedding of $\ss$ into $\gr$ is easy.
\smallskip

(2)\qua $A\inc B$ implies $\ba\inc\bar{B}$, so if $\ba\not=\bar{B}$ then 
the restriction of $x_{\bar{B}}$ to $\bar A$ is constant, so 
the restriction of $\bx_B$ to $A$ is also constant.

If $\ba=\bar{B}$, then $\bx_A$ and $\bx_B$ are both restrictions of $x_\ba$, so
$\bx_A$ is the restriction of $\bx_B$ to $A$.
\smallskip

(3)\qua The existence and uniqueness of such an $S$ come from the following
construction of $S$ as a function of an $x\in\gr$.

The elements of $S$ can be enumerated recursively, starting
with its first element $A=V$. 
For any $A\in S$, consider the partition of $A$ defined by $x_A$
(i.e. the set of preimages of the singletons). Then, the greatest
elements of $S$ strictly included in $A$ must be precisely the 
connected components, with cardinality greater
than one, of the sets in this partition. \endproof

\subsection{Description of the corners of $\gr$}\label{corn}

Now we identify $\ss$ with a subset of $\h$ thanks to Lemma \ref{plo}.
A direct calculation shows that dim 
$\ss=\dim C(G)-\#S'$. We are going to see the following

\begin{prop} $\ss$ has
a neighbourhood\footnote{A neighbourhood of a subset $P$ of a 
topological space is a neighbourhood of every point of $P$; 
it does not necessarily contain the closure of $P$.} 
in $\gr$ which is diffeomorphic to a neighbourhood of 
$\ss\times\{0\}^{S'}$ in $\ss\times[0, \1[^{S'}$. 
Thus, $\h$ is a manifold with corners.\end{prop}

Intuitively, in the family of parameters $(u_A)_{A\in S'}$, each $u_A$ will
measure the relative size of $A$ and $\ha$.

This diffeomorphism
is not completely canonical: to define it we have to choose for each $A\in S$
an identification of $\co(A/S)$ 
with a smooth section of representatives in $(\ra)^{A/S}$ where
$\ra$ is just a copy of $\rt$ marked with the label $A$. 
For instance, suppose that $V$ is totally ordered and that
for all $A\in R$ we fix the smallest vertex in $A$ at the origin of $\ra$; as for the
dilations, just fix any choice of smooth normalisation.
\ppar

{\bf Construction of a smooth map from a neighbourhood of 
$\ss\times\{0\}^{S'}$ in $\ss\times[0, \1[^{S'}$ to a neighbourhood of
$\ss$ in $\gr$}

Take a family $u=(u_A)_{A\in S'}$ of positive real numbers and a family
$f\in \ss$, that is, $f=(f_A)_{A\in S}$ and $f_A\in \ra^{A/S}$. Then, define
the net $(\phi_A)$ of correspondences from each space $\ra$ for 
$A\in S'$ to the space $\rb$ by
$$\eqalign{\phi_A:\ra&\longrightarrow\rb\cr x&\longmapsto f_\ha(A)+u_Ax.\cr}$$
We construct an element $g$ of $\h$ in the following way:

For all $B\in R$, we define the $C^B$ part of $g$ as a map 
from $B$ to $R_\bb$: first take for all $v\in B$ its image by
$f_{\overline{\{v\}}}$ in $R_{\overline{\{v\}}}$, then compose all the
maps $\phi_A$ above for the $A\in S$ such that
$\overline{\{v\}}\inc A\mathrel{\raise 1pt
\vtop{\hbox{$\subset$}\kern -6.7pt\hbox
{$\,\scriptscriptstyle\not=$}
}}\bb$ 
to obtain an element of $R_\bb$.

This map is not constant on $B$ for small enough values of $u$ because for
$u=0$ it is equal to $f_\bb$ which is 
not constant on $B$, so it gives an element of $C^B$. 

\def\U{{\cal U}}
This provides a well-defined smooth map from a 
neighbourhood $\U$ of $\ss\times\{0\}^{S'}$ in $\ss\times[0, \1[^{S'}$, to a
subset of $\h$, and its restriction to $\ss\times\{0\}$ is the imbedding
defined in Lemma \ref{plo}.

It maps the interior of $\U$ to the stratum $\co(G)$ of $\h$, because when 
$\qs A\in S'$, $u_A>0$, these operations are equivalent to first 
constructing the projection of $g$ in $C^V$, then restricting it to each
element of $R$, 
because all the maps $\phi_A$ are dilations-translations of
$\rt$. So it maps $\U$ to $\gr$.

\begin{remark} This diffeomorphism maps any $(f, u)\in \U$ into 
the stratum $\co(G,S_u)$ where $S_u=\{A\in S\,|\,A=V\hbox{ or }u_A=0\}$.
\end{remark}

{\bf Proof that the map above has a smooth inverse}

We have to check that the $u_A$ ($A\in S'$) and the $f_A$ ($A\in S$)
can be expressed as smooth functions of a $g$ in a suitable neighbourhood of $\ss$ in $\h$.

$f_A\in\co(A/S)$ is the modification of $g_A$ obtained by mapping 
each $A'\in S$ such that $\widehat{A'}=A$ to 
$g_A(a')$ where $a'$ is the smallest element of $A'$.
This $f_A$ maps the two ends of every edge of $A/S$ to two 
different points because it is in a suitable neighbourhood of $C(G,S)$ in $\h$.

$u_A$ is the suitably defined size of $A$ in $g_{\ha}$.
 \endproof

The description of $\gr$ is now finished and we can proceed with the 
compactification of $\c$.

\subsection{Construction of the compactified space $\cb$ of $\c$}

{\bf Definition of $\cb$}

Starting with a diagram $\g\in\dm$, define the graph
$G$ with the same set $V$ of vertices as $\g$, and define the set 
of edges of $G$ to be
the set of edges of $\g$, plus all pairs of univalent vertices.

First, canonically imbed the space $\c$ into the space
$$\hh=\h\times M^U$$
and let the compactified space $\cb$ of $\c$ be its closure in $\hh$. 

Now we are going to describe the corners of this space. For convenience,
we suppose that $M$ is oriented and $\del M=\0$. The case of 
$\del M\not=\0$ will be seen at the end of this subsection.

The reason why we took all pairs of univalent vertices as edges of $G$ and
not only the consecutive ones is to have $U$ connected even if $M$ is not, 
so that for any $(f, f')\in \cb\inc\h\times M^U$ such that $f'$ is not the constant
map to $M$, the restriction of the map $f_{\bar U}$ to $U$ is not constant 
and therefore is a picture of $f'$.
(For a knot, it makes no difference.)
\ppar

{\noindent \bf List of strata}

For commodity, let us fix an orientation on the manifold $M$.

A stratum of $\c$ is labelled by the following data of 
$(S,P,(\leq)_{P})$:

\b\qua A set $S\inc R$ as above (Notation \ref{nots})

\b\qua A partition $P$ of $U$ (that will be the set of preimages of $f'$), 
which satisfies one of the relations:

(i) $P=\{U\}$

(ii) $P=$ the partition of $U$ induced by $\bu/S$.

\b\qua A total order on each element of $P$.

\noindent such that: for any $A\in P$ and $B\in S$, $A\cap B$ with its total 
order is a set of consecutive elements for $\s$ respecting the orientation
of $M$ (Thus, $A$ is mapped to a 
single connected component of $M$, and there is only one possibility for 
this total order if $A$ is not the whole preimage of a connected component 
of $M$).\ppar

{\noindent \bf Description of strata}

The stratum with label $(S,P, (\leq)_{P})$
is the set of elements $(f, f')\in \h\times M^U$ such that:

\b\qua $f'$ belongs to the closure of $\s$ and maps two univalent vertices
to the same element of $M$ if and only if they belong to the same element of
$P$.

\b\qua $f\in\ss$

\b\qua Each ordered pair $(x,y)$ of elements of $U$ which
belong to the same element $U_1$ of $P$ and are consecutive for $\leq_{U_1}$
is mapped by $f_{\{x,y\}}$ to the 
tangent vector to $L$ at their common image. Or equivalently, $\qs
A\in S$, if $A\cap U$ is contained in one element $U_1$
of $P$, then $f_A$ maps $A\cap U$ to the straight line $\D$ with direction
the tangent vector $s$ to $L\circ f'(U_1)$, preserving the order 
in the large sense.

\b\qua In case (ii), $f_\bu$ must coincide on $U$ with $L\circ f'$. Thus,
we identify $R_\bu$ to the image space $\rt$ of $L$, because 
$L\circ f'$ and $f_\bu$ are not constant on $U$.
\ppar

Note that in case (ii), the vertices in $V-\bu$ are those which
escape to infinity, or for which any path from them to a univalent vertex
passes through trivalent vertices which escape to infinity 
(in case (i) there can be an undetermination).

We suppose that the univalent vertices are smaller so that they have priority
to be at the origin of the spaces $R_A$.\ppar

{\noindent \bf The neighbourhood of a stratum}

This is again
a corner, with a family of independent positive
parameters which are the same as before (indexed by $S'$), plus one 
more in the case (i), which will be denoted by $u_0$. 

To be quick, we shall only define the elements of $\c$ which correspond to
a family of strictly positive coordinates $(u_A)_{A\in S'}$, and 
possibly $u_0$. 

This will be made in two steps: first an approximate definition, then
a correction. 

In the first step, we start with an element $(f,f')$ of the considered stratum
with $f\in\ss$, and a family $(u_A)$ of small enough strictly positive
real numbers.  The construction of Subsection \ref{corn} gives an identification of 
all the spaces $\ra$ together (because $u_A>0$ for all $A$), 
and thus an element of $\co(G)$. There remains to identify one of the 
spaces (namely $R_\bu$) to the image space $\rt$ of $L$.
In case (ii), this identification is already done; 
whereas in case (i), it will be $$\eqalign{R_\bu&\lra\rt\cr x&\longmapsto L\circ f'(U) +u_0x.}$$
In the second step, we have to choose for each $\al$ in the image of
$L\circ f'$ a diffeomorphism $\ph_\al$ between two neighbourhoods of $\al$ 
which verifies the two conditions: first, it
approximates the identity near $\al$ up to the first order; second, it maps
the tangent line to the curve $L(M)$ at $\al$, to the curve $L(M)$ 
itself. This system
of diffeomorphisms must depend smoothly on $f$.

Then for every $\al$, correct the map from $V$ to $\rt$ of the first step 
by composing it with the map $\ph_\al$ for the vertices in 
${f_\bu}^{-1}(\al)$.

It can easily be seen that the resulting map extends as a smooth map
on a neighbourhood of the stratum.\ppar

{\noindent \bf Case when $\del M\not=\0$}

To the data $(S,P, (\leq)_{P})$ we must add the datum of the set 
$f'(U)\cap\del M$. As for the parametrisation of the corners, each $x\in 
f'(U)\cap\del M$ gives the parameter dist$(x,f'_1(U))$ for the $(f_1,f'_1)$ in
the neighbourhood of the stratum.
 \endproof
\ppar

{\noindent \bf End of proof of Proposition \ref{compac}} 

Now that we have compactified
$\c$ into a manifold with corners, we just have to check that $\Psi$
is smooth on $\cb$. But $\Psi$ is just the map defined by the
canonical projections of $\h$ on the $C^e\approx \sd $ for $e\in E$, 
since $E\inc R$. \endproof

\subsection{List of faces}

The codimension 1 strata of $\cb$ as constructed in the previous subsection
can be classified into six types. In the first five 
types we have $f'(U)\cap\del M=\0$.

(a)\qua In case (i) there is the coordinate $u_0$, so $S$ must be reduced to
$\{V\}$.

In case (ii) with $S=\{V,A\}$, let us distinguish four types:

(b)\qua $U\inc A$ (this corresponds to the case when some vertices go
to infinity),

(c)\qua $\#A=2$, $U\not\inc A$ and $E'_A\not=\0$; this is divided into two subtypes:

\quad (c1)\qua $A\inc T$ (thus $A\in E$),

\quad (c2)\qua $A\not\inc T$ (thus, either $A\subset U$ and $A\not\in E$, or $A\not\inc U$ and $A\in E$).

(d)\qua $\#A>2$, $U\not\inc A$ and $E'_A\not=\0$,

(a')\qua $E'_A=\0$ (thus, $A\cap U\not= \0$ and $U\not\inc A$).

Finally, 

(e)\qua Case (ii) with $S=\{V\}$ and $\#(f'(U)\cap\del M)=1$. 

\begin{defin} The faces of types (b), (c) and (d) will be called 
{\em ordinary faces}.

The faces of types (a) and (a') will be called {\em anomalous faces}\footnote{They are the ``anomalous faces'' of Bott and Taubes.}.

The faces of type (e) will be called {\em extra faces}; they do not exist in the case of links which is the main interest of this article.

We say that a face $F$ is {\it degenerated} if $\codim\Psi(F)>1$.
\end{defin}

\begin{nota} For $A\in R$ with $U\not\inc A$, let $F(\g,A)$ denote the face 
of $\cb$ defined by $S=\{V,A\}$.
\end{nota}

In the next section we shall prove the following proposition about
some types of faces which are degenerated.  Its aim is to eliminate
these diagrams for the next sections, as only the non-degenerated faces
will need to be taken into account.
For the same reason, we shall restrict the study to the configuration
spaces of subprincipal diagrams, as the faces of other spaces are all degenerated
thanks to Proposition \ref{34}.

\begin{prop}\label{nondeg} The only faces of $\cb$ which can be non-degenerated
are:

The anomalous faces in the connected case (type (a) faces with $\g$ connected, 
and type (a') faces $F(\g,A)$ with $A$ connected).

The faces of types (c) and (e).

The faces $F (\g,A)$ of type (d) which satisfy the following conditions: each 
edge in $E'_A$ meets at least two edges in $E_A$, and the number $\#E'_A$
is $1$ or $2$ if $A\not\inc T$, and $3$ or $4$ if $A\inc T$.

\end{prop}

\section{Proofs for the codimensions}\label{pr34}
\setcounter{equation}{0}
\subsection{Necessary conditions in Proposition \ref{34}}
Let $A\inc T,\ \#A>1$. 
Note that in the definition of principal or subprincipal diagrams, 
we can equivalently
restrict the conditions to the sets $A$ which are connected: replacing
a disconnected $A$ with one of its connected components does the trick.

Consider the following commutative diagram:
$$\diagram{\c&\hfl{\Psi}&(\sd )^{E}\cr\vf&&\vf\cr C^A&\hfl{}&(\sd )^{E_A}}$$
According to (\ref{count}) we have $2\#E_A+\#E'_A=3\#A$. Thus,
$$\dim C^A=\dim(\sd )^{E_A}+\#E'_A-4.$$
So, if $\Psi(\c)$ has codimension 0, then $\dim C^A\geq
\dim(\sd )^{E_A}$ for all $A$, thus $\g$ is principal.

This way, we can see that $\codim \Psi(\c)\leq 1$ implies Condition (1) 
in the definition of a subprincipal diagram.
Similarly, considering for any disjoint $A,A'\inc T,\ \#A>1,\ \#A'>1$, 
the natural map from $\c$ to $(\sd )^{E_A}\times (\sd )^{E_{A'}}$,
we can see that it also implies Condition (2).
Therefore, all diagrams such that $\codim \Psi(\c)\leq 1$ are subprincipal.
 \endproof

\def\ccl{\overline{\cal F}}

\subsection{A general lemma on codimensions}
\label{dime}

In the proofs of the propositions, we shall use several times the following
argument to minorate the codimension of the images by $\Psi$ of certain
subsets of $\h$. It will especially apply 
to the properties of the image of a stratum of 
$\cb$ in $C^V\times M^U$
by the canonical projection (values of $(f_V, f')$). It will play the role of
a lemma but its expression is long and its proof trivial.

\def \pp{{\cal P}}
Let $\g$ be a connected, subprincipal diagram.

Let $\pp\inc R\backslash\{V\}$ be a set of pairwise disjoint sets of vertices.

Let $\g'=\g/\pp$ be the quotient diagram, made of 
\smallskip

\b\qua a set $V'=U'\sqcup T'$ of
vertices such that:

(1)\qua $V'$ is the quotient of $V$ by the relation whose equivalent classes are the elements of $\pp$
and singletons, 

(2)\qua $U'$ is the image of $U$ by the canonical projection from $V$ to $V'$.
\smallskip

\b\qua a set of edges $E'=\{e\in E\,|\,\qs A\in\pp,\ e\notin E_A\}$.
\smallskip

We can easily see that $\qs x\in U', (x\in \pp$ or $(x\in U \hbox{ and is univalent in }\g'))$.

The same for $T'$: $\qs x\in T', (x\in \pp$ or $(x\in T \hbox{ and is trivalent in }\g'))$.

Let $U''=U'\cap\pp$ and $T''=T'\cap\pp$.

Let $L$ be a curve in $\rt$, and let $C(\g',L)$ be
some configuration space of $\g'$ where every element of $U'$ runs along $L$.

For any $x\in V'$, let $n_x$ be the number of edges in $E'$ which contain
$x$. Then, $$\dim(\sd)^{E'}-\dim(C(\g',L))=\sum_{x\in T''}(n_x-3)+\sum_{x\in U''}(n_x-1).$$
But (as $\g$ is connected and subprincipal), we have
$n_x-3\geq 0$ for any $x\in T''$, and $n_x-1\geq 0$ for any $x\in U''$, thus 
$\dim(\sd)^{E'}\geq\dim(C(\g',L))$.

Therefore, the codimension $c$ of the image of the canonical map (restricted $\Psi$) from
$C(\g',L)$ to $(\sd)^{E'}$ has the following properties:

(1)\qua If $c=0$ then the group of translations or dilations letting $L$ invariant has dimension 0, 
$n_x=3$ for every $x\in T''$ and $n_x=1$ for every $x\in U''$. Thus, as $\g$ is subprincipal,
every element of $U''$ contains at least two elements of $U$, and if $\g$ is principal then $T''=\0$.

(2)\qua If $c=1$ then the above consequences of $c=0$ can suffer only one exception with the range of one unit.

In particular, if $L$ is contained in a straight line, then $c>1$.

\subsection{Proof of the sufficient condition in Proposition \ref{34}}

\def\sdn{(\sd)^E}

We shall prove now that if $\g$ is principal then there exists a link
$L$ such that 
$\codim \Psi(\c)=0$. (This is just an elegant result which is not useful for 
the rest of the present article).

So, we will study the existence of preimages for generic points in $(\sd)^E$.

\def\ce{{\cal E}_s}

For any $s=(s_e)_{e\in E}
\in \sdn$, let $$\ce=\prod_{e\in E} \rt/\R s_e.$$
Let $\Psi_s$ be the linear 
map from $(\rt)^V$ to $\ce$ defined by: for each edge $e=(v,v')\in E\inc V^2$, 
and each $f\in (\rt)^V$, $\Psi_s(f)_e \in \rt/\R s_e$ is the image of 
$f(v')-f(v)$ by the canonical linear projection.

Now, the configuration space $\c=\c(L)$ of a link $L$ is a submanifold of
$(\rt)^V$ with the same dimension as $\ce$, 
and we shall consider the restriction of $\Psi_s$ to it. 
Then, the set of preimages of $0$ by $\Psi_s$ in 
$\c$ corresponds bijectively to the disjoint union of sets of preimages by $\Psi$ in $\c$ of all elements of the form $(\pm s_e)_{e\in E}$.

The condition ``$\codim \Psi(\c)=0$'' can be reformulated in the form
``There exists a configuration $x\in\c$ at which the differential map d$\Psi(x)$ 
is bijective''. 
Let us analyse this differential map. 

\def\bcl{\overline{C_{\g'}}(L')}
The tangent space of $\c$ at $x$ is the topological closure 
$\bcl\inc (\rt)^V$ of the configuration space $C_{\g'}(L')$ where $L'$ 
is the family of the tangent straight lines to the link $L$ 
at the positions of univalent vertices in the configuration $x$, 
and $\g'$ is obtained from $\g$ by replacing its support with $U\times\R$.

Then, d$\Psi(x)$ is the map $\Psi_s$ defined on $\bcl$ where $s=\Psi(x)$.

What we have to prove is, for any principal diagram $\g$, the existence 
of some $s,x$ and $L$ which verify these conditions. First, let us fix
a family $L'$ of pairwise disjoint straight lines indexed by $U$.

\begin{lemma}
For a generic $s\in \sdn$, the map 
$\Psi_s$ is bijective from the affine space $\bcl$ to $\ce$.
\end{lemma}

\proof
If it was not, as $\bcl$ and $\ce$ have the same dimension, then the kernel 
of the linear part of $\Psi_s$ would be nonzero, that is,
there would be a nonzero element $x\in \overline{C_{\g'}}(TL')$ 
(where $TL'$ is the family of parallel lines to the components of 
$L'$ at the origin), in Ker $\Psi_s$.
But this is not possible for a generic $s$, thanks to Subsection \ref{dime}, 
because $TL'$ is invariant by dilations with center the origin.
 \endproof

Now, $0$ has one preimage $x$ by $\Psi_s$ in $\bcl$ for a generic $s$. This
$x$ is not collapsed to one single point of $\rt$ because the lines in $L'$
are disjoint.
Let us apply once more the result of Subsection \ref{dime}: as we suppose
that $\g$ is principal, we find $T''=\0$,
and also $U''=\0$ because all univalent vertices run over pairwise 
disjoint lines.

Therefore, we have $x\in C_{\g'}(L')$.
Finally, we can find a link $L$ which is tangent to every line of $L'$
at the position of every univalent vertex of $x$ respecting the datum $\sigma$
of $\g$, so that at the same geometric configuration $x$, we have
d$\Psi(x)\approx \Psi_s$ bijective from $T_x\c=C_{\g'}(L')$ to $T_s(\sd)^E
\approx{\cal E}_s$ where $s=\Psi(x)$.
 \endproof

\subsection{Proofs for the degenerated faces}

\begin{remark} \label{discon}
Suppose there exists a subset $A$ of $V$ such that $E'_A=\0$. Then
consider the two subdiagrams $\g_1$ and $\g_2$ whose sets of vertices are
$A$ and $V-A$ respectively. 
Then there is a canonical diffeomorphism from $\c$ to an open subset of
$C_{\g_1}\times C_{\g_2}$, which is 
delimited by a finite union of subspaces. 
\end{remark}

This gives a diffeomorphism between a type (a') face, 
and an open subset of the product of a type (a) face of $C_{\g_1}$ with
the space $C_{\g_2}$. So it shows that the types (a) and (a') are fundamentally
the same.

It also has the following consequence:

\begin{lemma}
A type (a) face with $\g$ not connected, or a type (a') face
with $A$ not connected, is degenerated.
\end{lemma}

\proof
When $\g$ is a disjoint union of two diagrams $\g_1$ and $\g_2$,
the natural map from $\cb$ to $\overline {C_{\g_1}}\times \overline {C_{\g_2}}$ 
maps a type (a) face of $\c$ to a product of type (a) faces of 
$\overline C_{\g_1}$ and $\overline {C_{\g_2}}$ , which has codimension 2. Thus, as $\Psi$ is 
factorised through this, the image of its restriction to this face
also has codimension at least two.

This concludes for the type (a), and the generalisation to (a') is immediate. \endproof

\begin{lemma}\label{btype}
The type (b) faces, and the type (d) faces $F (\g,A)$ such that  
($A\not\inc T$ and $\#E'_A>2$) or ($A\inc T$ and $\#E'_A>4$), 
are degenerated.
\end{lemma}

These are direct consequences of Subsection \ref{dime}. For the type (b),
the image of the univalent vertices is a vertex at least trivalent (because
the diagram is subprincipal) which must stay at the origin, so that
the codimension is at least 3. \endproof

\begin{lemma}\label{facedil}
Any type (d) face such that the vertex $v$ in $A$ of some
edge in $E'_A$ is either in $U$ or not bivalent in $A$, 
is degenerated.
\end{lemma}

Indeed, in such a face, each preimage of a point by $\Psi$ is at least 
one-dimensional because in the $C^A$ part of this preimage, this 
vertex $v$ can move freely in one direction while the others remain fixed.
This direction is the tangent direction of $L$ at $f'(v)$ if $v\in U$,
or the one of its edge in $E_A$ if $v$ is univalent in $A$.
This move is not eliminated by dilations because $\#A>2$.
 \endproof

This ends the proof of Proposition \ref{nondeg}.

\section{Labelled diagrams and orientations}
\setcounter{equation}{0}

\subsection{Labelled diagrams}\label{labd}

To formulate the degree $n$ part of $Z^0$ in terms of the degree of the map
$\Psi$ from the glued configuration spaces to a product of spheres, we shall
need to fix this product of spheres. Thus, all our edges will be
labelled and oriented. We shall call that ``labellings"
of diagrams.
We are going to reformulate the definition of the diagrams so that they
will be intrinsically labelled.

Let us fix the degree $n$ of the diagrams. Their number of edges 
varies under STU but we must fix the number $N$ of their labels.
So, we shall only work with the diagrams $\g$ which have a number of edges 
$e_\g\leq N$, and thus a number of univalent vertices $u_\g\geq k=3n-N$, 
and there may be some unused labels which will be called ``absent edges".

\begin{nota}If $A$ is a set of sets, let $\cup A$ denote the union of 
elements of $A$.

Let $n\geq 1,\ k\leq2n$ be fixed integers and 
$N=3n-k$.

Let $E=\{e_1,\cdots, e_N\}$ be the set of labelled ``edges" 
$e_i=\{2i-1,2i\}$.
The set $\dem=\cup E=\{1,\cdots, 2N\}$ is called the set of half-edges.
\end{nota}

\begin{defin}\label{diag}
Let $\dl$ denote the set of 4-tuples 
$\g=(U,T,E^v,\s)$,
called {\it labelled diagrams}, 
which respect the following conditions:

\b\qua $E^v\inc E$ is the set of ``visible edges", it labels the edges of 
the diagram. We let $E^a=E-E^v$ denote the set of ``absent edges".

\b\qua The set $V=U\cup T$ of all vertices of $\g$ is a partition of $\cup E^v\inc
\dem$.

\b\qua $U$ is the set of univalent vertices: it is a set of singletons, 
and $\s$ is an isotopy class of injections of $U$ into the interior of $M$.

\b\qua $T$ is the set of trivalent vertices, which are subsets of $\dem$ with 
three elements belonging to three different edges.

\b\qua $\#V=2n$ (which is equivalent to $u_\g=\#U=\#E^a+k$).

\b\qua $\g$ is subprincipal (see Definition \ref{3con}).
\end{defin}

The identity is the only isomorphism of diagrams which preserves 
the labelling, so we will take the same notion of an isomorphism as before,
which refers to the underlying unlabelled diagrams:

\begin{defin} An {\em isomorphism} (or {\it change of labelling})
between two labelled diagrams 
$\g=(U_\g,T_\g,E^v_\g,\s_\g)$ and $\g'=\-(U_{\g'},T_{\g'},E_{\g'}^v,\s_{\g'})$ 
 is a bijection from $\cup E^v_\g$ to $\cup E_{\g'}^v$ which carries
$E^v_\g$ to $E^v_{\g'}$, $T_\g$ to $T_{\g'}$, $U_\g$ to $U_{\g'}$ and 
$\s_\g$ to $\s_{\g'}$.
\end{defin}

A configuration for a labelled diagram induces a map from $\cup E^v$ 
to $\rt$, so the canonical map $\Psi$ can be written: 
$$\eqalign{\Psi:C_\g&\longrightarrow (\sd )^{E^v}\cr 
\f&\longmapsto
\left({\f(2i)-\f(2i-1)\over|\f(2i)-\f(2i-1)|}\right)_{e_i\in E^v}
\cr}$$
Then, we shall denote by $\Psi'$ (but also sometimes by $\Psi$ 
by abuse of notation; fortunately, all ambiguities will be ineffective to
the results)
the map from $\cc=\cb \times (\sd )^{E^a}$
to $(\sd )^E$ defined by the product of the previous $\Psi$ with the identity
on$(\sd )^{E^a}$.

\subsection{Orientation of the configuration spaces}\label{ori}

We shall now define an orientation of the 
space $\cc$ depending on the orientation $o$ of $\g$. It will be done by
labelling all components of a natural local coordinate system of $\cc$ 
by the elements of $\dem$. These components are:

\b\qua One coordinate in $M$ for each univalent vertex, respecting the local
orientation of $M$ given by the orientation of this vertex.

\b\qua The three standard components in $\rt$ for each trivalent vertex.

\b\qua And for each absent edge we have the two components of a local 
coordinate system in $\sd $ with the standard orientation (i.e. where a 
basis $(x_1,x_2)$ at a point of $\sd$ is positively oriented when the 
vector product $x_1\wedge x_2\in\rt$ points outwards $\sd$).
\smallskip

The coordinate in $M$ for a univalent vertex $\{x\}$
is labelled by $x$. The two components for each $e_i\in E^a$
are respectively labelled by $2i-1$ and $2i$.

As for the trivalent vertices, let us represent
$o$ by a family of bijections $(b_t)_{t\in T}$ from each trivalent vertex $t$
to the set $\{1,2,3\}$.
Then, for all $x\in t$, give the label $x$ to the $b_t(x)$-th component of 
$\f(t)$ in $\rt$.

Thus we have an orientation of $\cc$, which only depends on the orientation
$o$ and the labelling of $\g$.
It is reversed when we change the orientation of one vertex or one edge;
and the isomorphism induced by a permutation of $E$, which preserves the 
orientations of edges and vertices, also preserves the orientation of $\cc$. 
So, this construction is equivalent to
first orienting $\c$ by the same method with edges $e_i$ 
($i=1,\cdots,\#E^v$), then taking the direct product with $(\sd)^{E^a}$.

For a labelled oriented diagram $(\g,o)$, we can now define the integral 
$$I_L(\g,o)=\int_\cc {\Psi'}^*\O=\int_\c \Psi^*\O_{E^v}$$
where $\O$ is the standard volume form on $(\sd )^E$ with total mass 1,
and $\O_{E^v}$ is the one on $(\sd )^{E^v}$.

Note that we have $I_L(\g,-o)=-I_L(\g,o)$ and that $I_L(\g,o)$ is independent 
of the labelling of $\g$: exchanging two edges preserves both 
orientations of $\cc$ and $(\sd )^E$, and reversing one edge changes both 
orientations (precisely, it changes both the orientation of
$\cc$ and the sign of the form ${\Psi'}^*\O$), because the antipodal 
map reverses the orientation of $\sd $.

So, the product $I_L(\g)[\g]$ is well-defined independently 
of $o$.

In all the following, we fix an arbitrary choice of orientations 
on the diagrams in $\dl$, 
with no consequence on the final results.

\subsection{Interpretation of $Z$ in degree theory}\label{degthe}
\def\C{C'}

Let $\poids $ be a map from $\dl$ to $\ank$
such that: 
\begin{equation}\label{somme}
\qs\g\ \hbox{principal},\qua \sum_{\g'\in\dl,\g'\sim\g}\poids(\g')={[\g]\over |\g|}
\end{equation}
The standard example will be 
\begin{equation}\label{coef}
\poids (\g)=\hbox{projection of }{(N-e_\g)!\over N!\;2^{e_\g}}[\g] \hbox{ in } \ank
\end{equation}
but we need not restrict to it now.

Let $\C$ be the following formal linear combination of spaces
\begin{equation}\label{space}
\C=\sum_{\g\in \dl}\poids (\g) \cc.\end{equation}
Here, $-\cc$ stands for $\cc$ with the opposite orientation.
Let $\xi$ be the almost everywhere defined map with rational values 
in $\ank$:
$$\eqalign{\xi:(\sd)^E&\lra\ank\cr x&\longmapsto 
\hbox{ count of }\Psi^{-1}(x)\hbox{ in }\C }
$$
According to the degree theory, this map can be extended 
into a locally constant map defined
outside the union of the $\Psi(\del\cc)$ for $\g\in\dl$. When a
map is defined on a manifold without boundary, the number of preimages
counted with signs is constant and called the {\it degree} of the map.
In this article, we shall use a kind of generalised notion of ``manifold 
without boundary'', with a notion of gluing of faces that we shall define
in the next section, such that
only the faces which are not glued are relevant in the notion of a boundary.
Now, we have $$Z_n^k=\int_{(\sd)^E} \xi\O.$$
So, $Z_n^k$ is generally not rational because the map $\xi$ is not 
constant on $(\sd)^E$.
The aim of the following sections will be to prove Propositions \ref{varia} and
\ref{mult} in the following way: we shall introduce other configuration spaces to 
complete $\C$ 
so that the corresponding $\xi$ becomes constant. 

Then, each of the expressions $Z^0$ of Proposition \ref{varia} and 
$Z_n^k$ of 
Proposition \ref{mult} will be equal to the degree of the map $\Psi$
defined on such a linear combination of configuration spaces which extends $\C$.

\begin{remark} We could generalise these constructions to other maps 
$\poids $ from $\dl$ to any $\Z$-module. 
This possibility will be discussed in the Appendix. We shall see that
in the case of a torsion-free module it gives no more invariants. The
only interest is that it can improve the rationality result by means of
another map $\poids$ which satisfies Formula 
(\ref{somme}) and optimised gluing conditions, but not
(\ref{coef}). 
However, the improvement expressed in Remark \ref{improve} is not based on this method. 
\end{remark}

\subsection{Yang's point of view}\label{yang}

Yang \cite{ya} introduced another interpretation of $Z$ in 
terms of degree theory,
which is sometimes useful 
to avoid some computations of coefficients. 
Let us keep here Formula (\ref{coef}) for $\poids$. 

\def\gg{{\cal G}}\def\rp{{\sym RP}}
The idea is to quotient the space $(\sd)^E$ by the group $\gg$ of 
permutations and reversals of the edges, to obtain the space 
$B_E=(\rp^2)^E/\Sigma_E$ where $\Sigma_E$ is the permutation group of 
$E$. The group $\gg$ also naturally acts on $\C$, so that $\gg$ 
commutes with $\Psi$ in these actions.
Thus, the map $\xi$ is invariant under the action of $\gg$ (for the same reason
as the fact that $I_L(\g)[\g]$ did not depend on orientations) and induces
a map defined on $B_E$.

To understand the quotient space $\C/\gg$, first note that the
group of automorphisms of a diagram $\g$ freely acts on $\c$.
Let $\c/|\g|$ denote the quotient space of $\c$ by this action.
It is almost everywhere identical to
the space of all ``configured diagrams'' isomorphic to $\g$, that is, 
unlabelled diagrams where the set $V$ of vertices is a subset of the 
space $\rt$ and the univalent vertices belong to the image of $L$.
Then $\C/\gg$ is topologically the union over all isomorphism 
classes of degree $n$ diagrams $\g$ with $u_\g\geq k$, of the product 
$\c/|\g| \times B_{E^a}$.

Now, the expression $Z_n^k$ can be interpreted as an integral on these quotients,
but not of a $2N$-form because some $I_L(\g)$ can contribute even when
the quotient space $\c/|\g|$ is not orientable (the simplest example is 
$\g=\th$). So we must replace the form $\Psi^*\O$ by the density
$(\Psi^*\O\cdot($local orientation of $\c))$. 

Let us drop the condition $u_\g\geq k$ (or equivalently, take $k=2$);
replace $\C$ by the space $C$ constructed the same way as $\C/\gg$
but without the absent edges (it is the union of the $\c/|\g|$ for all 
degree $n$ diagrams $\g$ up to isomorphism). 
Thus it is made of parts with different dimensions.

So, we can say that $$Z_n=\int_C \left(\Psi^*\O\cdot(\hbox{orientation})\cdot[\g]\right)_{local}.$$
The whole integral $Z$ itself is defined in the same way with diagrams
of all degrees.

Note that 
the rational coefficient disappears in the above integral formula. 
The advantage of this presentation is that it generally trivialises the problem
of the handling of the rational coefficients of Formula (\ref{coef}). The bad point is that 
it gives no means to handle the sign of the density involved in
the space $C$. 

\section{Ordinary gluings}\label{recol}
\setcounter{equation}{0}

\subsection{Definition of a gluing}

\def\glu{$\Psi$-isomor\-phism\ }
\def\glus{$\Psi$-isomorphisms\ }
In this subsection, we explain in all generality (for any map $\poids$
and without restricting to the set of configuration spaces we already defined)
the notion of a gluing of a set of faces:

Each face $F$ of a configuration space $\cc$ contributes to a 
variation of $\poids(\g)$ of the map $\xi$ while going across $\Psi(F)$.
Here, $\Psi(F)$ is counted with possible multiplicities.
We are going to define {\it gluings} for certain sets of faces, 
that is, objects whose existence ensure that the contributions of these faces
to the discontinuity of $\xi$ at their images
cancel. Note that the degenerated faces do not contribute to the discontinuity of
$\xi$, thus no gluing is necessary for them.

\def\M{{\cal M}}
\begin{defin} A {\em $\Psi$-isomorphism} between two manifolds $\M_1$ and $\M_2$ 
with maps $\Psi$ from each
of them to $(\sd)^N$, is
a diffeomorphism $\ph$ from $\M_1$ to $\M_2$ such that $\Psi\circ \ph=\Psi$ 
on $\M_1$. A $\Psi$-isomorphism is said to be {\em positive} if it preserves their
orientations.
\end{defin}

\begin{defin} A {\it self-gluing} of a face is a negative $\Psi$-automorphism
of this face.
\end{defin}

The existence of a self-gluing of a face clearly implies that this face does
not contribute to the discontinuity of $\xi$. But self-gluings are quite
exceptional and we shall use them only to eliminate some pathological cases
of the other gluings.

In this section, we define gluings in the following way 
(but for the anomalous faces in the next section, 
it will be necessary to cut some faces into pieces).

\begin{defin} A {\it gluing} of a set of faces $F (\g,A)$ with coefficients 
$\poids (\g)$ is a partition of this set such that in each class, all 
the faces  are $\Psi$-isomorphic and verify the {\it gluing condition}
$$\sum_{(\g,A)\hbox{ in the class}} \pm \poids (\g)=0$$
where the signs in front of $\poids (\g)$ and
$\poids (\g')$ for two given pairs $(\g,A)$ and $(\g',A')$ will be equal iff
each $\Psi$-isomorphism between $F (\g,A)$ and $F (\g',A')$ 
preserves their orientations inherited from the orientations
of the spaces $\cc$ and $C'_{\g'}$. 
\end{defin}

Note that if there are both a positive and a negative 
$\Psi$-isomorphism between two faces, then each face in the class is 
glued to itself, so the condition may be deleted. Otherwise, the condition for
the signs in the above formula is consistent.

\subsection{Gluing results for ordinary faces}

The aim of this section is to prove the following proposition. 

\begin{prop} \label{glunano}
For the map $\poids$ defined by Formula (\ref{coef}),
there exists a gluing of the set of ordinary, non-degenerated
faces of the space $C'$ defined by (\ref{space}).
\end{prop}

\proof We shall prove the following three lemmas
in Subsection \ref{proglu}. They apply to the gluing of a combination $C'$ 
as above but that can come from any map $\poids$ defined on the set 
of diagrams $\dl$ with chosen orientations, and seen as extended to 
the set of the same diagrams with all orientations, 
according to the AS relation.

\begin{lemma}\label{rihx} The gluings of the type (c1) faces are ensured
by the IHX' relations defined as follows: they are all possible relations
of the following form between six subprincipal labelled diagrams that are 
identical outside the drawn part:
$$\poids \epfb{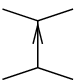}+\poids \epfb{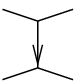}=\poids \epfb{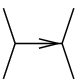}+\poids \epfb{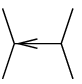}-\poids \epfb{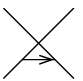}-\poids \epfb{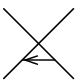}$$
\end{lemma}

\begin{lemma}\label{rstu} The gluings of the type (c2) faces are ensured
by the following STU' relation between labelled diagrams that are 
identical outside the drawn part and such that all three unlabelled diagrams 
of the corresponding STU relation are subprincipal:
$$\poids \epfb{u}-\poids \epfb{s}
=\sum_{e\in E^a}\left( \poids \epfb{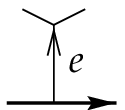}+\poids \epfb{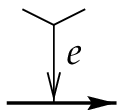}\right)$$
where $E^a$ is the set of absent edges of a diagram in the left-hand side 
of the relation (it can be empty).
\end{lemma}

If this $E^a$ is empty, this relation is verified in the case of Formula 
(\ref{coef}) because
it is deduced from the STU relation and the fact that the missing
diagram was forced to cancel in $\ank$.

\begin{lemma}\label{autorel} There is a gluing of the type (d) faces 
if $\poids $ is a constant on each isomorphism class of diagrams.
\end{lemma}

The relations stated in the first two lemmas above will be fully used in
the proofs to define the gluings, but the proof of the last
lemma will not make full use of its assumption (that $\poids $ is 
constant on each isomorphism class) but only of a consequence of it, whose
expression would be complicated to write and useless here. This possibility
to weaken the condition can be used to improve the result on the 
denominators of the integral (see Appendix
\ref{r2}).

Then, we finish the proof of Proposition \ref{glunano} by the remark that
the map $\poids$ defined by (\ref{coef}) satisfies the conditions of these
lemmas.

\subsection{Gluing methods and orientations}

\let\eps=\varepsilon
\def\sign{\mathop{\rm Sgn}}

Now we are going to consider a general case of \glu between two faces
$F (\g,A)$ and $F (\g',A')$ of a same type ((c) or (d)). The problem will be
to determine its sign $\eps$. 

We suppose that $\g/A=\g'/A'$ (as labelled diagrams), and that the \glu
preserves the projection onto $C_{\g/A}$. 
So, the only transformations concern the graphs $A$ and $A'$ and the 
spaces $\co(A)$ and $\co(A')$ (Notation \ref{notgra}), and if 
$A\cap U$ and $A'\cap U'$ are nonempty, 
they are mapped to the same line $\D$ in the two identified configurations.

Moreover, let us suppose that $E_A=E_{A'}$ (it will be the case for the 
gluings we shall consider except one of the two kinds of gluings for the 
STU' relation).
This implies that $A$ and $A'$ have the same respective
numbers of univalent and trivalent vertices.
\ppar

\def\L{{\cal L}}
For convenience, we are going to label all vertices in $A$ independently
of the labelling of the edges, starting with
the univalent vertices (this means the elements of $U$), if they exist.
We denote them by $x_i$ for $i\in X$, $X=\{0,\cdots,\#A-1\}$. 
Thus, $x_0\in U$ if $A\cap U\not=\0$. By convention,
we put the vertex $x_0$ at the origin of the space $R_A$. 
So $C^A$, which contains $\co(A)$, is identified
with the space of half lines in the vector space $(\rt)^{X-\{0\}}$. 
We do the same for $\g'$ (so the univalent vertices in $A$ and $A'$
are labelled by the same elements of $X$). 

The neighbourhoods of the faces are locally identified (on open dense
subsets, in natural ways with a definite limit near the interior of the face) 
with the fiber bundle
$\ee$ with base $C_{\g/A}$ and fiber $\ce$, where $\ce$ is the space defined by:
$\ce=(\rt)^{X-\{0\}}$ if $A\cap U=\0$; and if $A\cap U\not=\0$,
$\ce$ is the subspace of $(\rt)^{X-\{0\}}$ consisting of the configurations 
which map $A\cap U$ to $\D$, where $s$ is defined by the base point in 
$C_{\g/A}$.

We shall suppose that the bijection between the faces is defined 
over each element of $C_{\g/A}$ by a linear transformation $\L$ of $\ce$.
If $A\cap U\not=\0$, consider the set of vertices $A\cap U$ 
as a subset of $X$: we shall suppose that $\L$ is the identity on the
corresponding copies of $l_s$, and that the orientations of these
vertices are preserved from $\g$ to $\g'$.

Then, the bijection between the faces smoothly extends to a bijection of
their neighbourhoods. As the orientation of each face of a space $\cb$ was 
defined as coming from the orientation of $\cb$, the sign $\eps$ of the 
\glu between $F (\g,A)$ and $F (\g',A')$ is given by the action on 
the neighbourhoods of the faces.

These neighbourhoods are both identified to the same space $\ee$, but 
the orientation of 
this space $\ee$ depends on the diagram.
Therefore, we have $$\eps=\sign(\det\L)\cdot\eps'$$ where $\eps'$ is
the sign of the change of orientation
of $\ee$ as orientations of the respective configurations spaces to which
$\ee$ was identified.

This $\eps'$ can be computed by comparing these
orientations of $\ee$ coming from the respective diagrams as defined 
in Subsection \ref{ori}.
The local coordinates, first labelled by
the vertices, are then labelled by the half-edges. 
The labelled graphs $A$ and $A'$ respectively define two bijections $b$ and $b'$ 
between
$$\{i\in X\,|\,x_i\in U\}\cup\{(i,1),(i,2),(i,3)\,|\,i\in X,x_i\in T\}$$
and the set $\cup A\inc \dem$ of half-edges which belong to a vertex in $A$.

Then, $\eps'$ is the signature of $b^{-1}b'$.

\subsection{Proofs for the gluings}\label{proglu}

{\noindent\bf Proof of Lemma \ref{rihx}}

Each term of an IHX' relation corresponds to a type (c1) face
$F(\g,e)$ where $\g$ is the labelled diagram of the term, and $e$ is
the edge involved in the relation. 

If some type (c1) face does not fit into such a relation, it is because 
one other diagram of the IHX relation is not subprincipal, 
so its faces are degenerated. Since all 
faces of this unwritten IHX' relation are $\Psi$-isomorphic, they are
all degenerated. This is why no relation is necessary for them.

So we have a partition of at least all non-degenerated type (c1) faces, 
and what we have to check is that this partition is a
gluing, where the gluing condition is the IHX' relation.

In the above constructions, we can always take $\L=$ Id, preserve the
graph structure of $A$ but just use permutations of the set $E'_A$. 
This provides $\Psi$-isomorphisms of faces with sign the
signature of this permutation (it must be considered
modulo the transpositions of two edges attached to the same vertex
of $A$, because these transpositions induce the identity).

The verification of the signs in the IHX' relation is now easy. \endproof

{\noindent\bf Proof of Lemma \ref{autorel}}

We consider now the non-degenerated type (d) faces. 

According to Proposition \ref{nondeg}, each edge in $E'_A$ meets
two edges (denoted as pairs of vertices) $e_1=\{x,y\}$ and
$e_2=\{x,z\}$ in $E_A$. Suppose $x$ has
a non-zero label. Let $\L$ be the map
which maps (the position of) $x$ to $y+z-x,$ and preserves the
other vertices.
This transformation has determinant $-1$.
$$\epsfbox{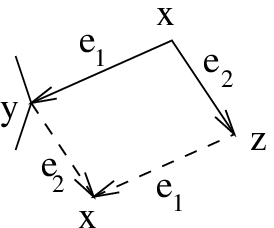}$$
Then apply the only change of structure of $A$ which fits in with this 
transformation: it is a mere relabelling, which exchanges $e_1$ and $e_2$ 
and reverses them: if $e_1$ goes from $x$ to $y$ in $A$, then it
goes from $z$ to $x$ in $A'$.
Suppose that these diagrams have the same orientation through this relabelling.
This defines an even permutation of $\cup A$, so this \glu is negative.

Now, here is a way to give a partition of these faces into pairs, which
is a gluing as soon as $\poids$ is constant on each isomorphism class of
diagrams with the same orientation.

Consider an (equivalence) class of pairs $(\g,A)$ (concerned here) 
which are isomorphic by isomorphisms whose restriction to $\g/A$
is the identity (thus $\g/A$ is fixed as a labelled diagram).
Then fix an arbitrary edge $e\in E'_A$. 
By this choice, the diagrams in this class are grouped in pairs
according to the above construction. \endproof

{\noindent\bf Proof of Lemma \ref{rstu}}

First let us check that for each pair $(\g, A)$ where $A$ is an edge
or a pair of consecutive univalent vertices of $\g$, the face $F (\g, A)$ can be 
identified with the space $C'_{\g/A}$:
indeed, let us distinguish two cases.

If $E_A=\0$, it comes from the fact that the configuration of $A$ 
on this face in the compactified configuration space, is 
determined by the one of $\g/A$ and the tangent map of $L$. 

If $E_A\not=\0$, we have
$A\not\inc U$. This allows us to put
the unique element of $E_A$ among the absent edges.

Therefore, all faces of a STU' relation (faces of the form $F(\g, A)$ 
where $\#A=2$ and $\g/A$ is identical to a fixed diagram with one 
bivalent vertex $\{A\}$ on $M$) are $\Psi$-isomorphic,
thus, they are all degenerated when one is. This is true for $k=2$ when
all faces of the STU relation are represented, thus, as the codimension does
not depend on $k$ (which determines the number of absent edges), all 
faces of an STU' relation are degenerated when one 
face of the STU relation with unlabelled diagrams is.

This proves that aevery non-degenerated type (c2) face 
appears in an STU' relation. 

Let us now study the signs of these $\Psi$-isomorphisms and check
that the STU' relations are the gluing conditions.
The sign of the \glu between both left-hand side diagrams comes from the
same method as in the IHX' gluing (here the space $\ce$ has dimension 1, 
$\det\L=1$ and $\eps'=-1$ because of the transposition). 

Otherwise, let us study the sign of one \glu between two faces corresponding to diagrams with a different number of univalent vertices.
This is the exceptional case of sign study, but we shall treat it in a similar
way to the other cases.

Consider the following diagrams, where each half-edge is marked with the
name of a component in a local coordinate system of the neighbourhood
of each face.
Here $(x,y)$ stand for the spherical coordinates, and $u$ stands
for the radial coordinate (outward normal vector to the sphere, thus 
inward normal vector to the face of the configuration space):
$$\epsfbox{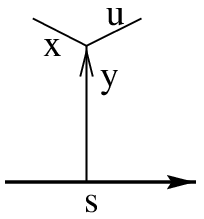}\hskip 3cm\epsfbox{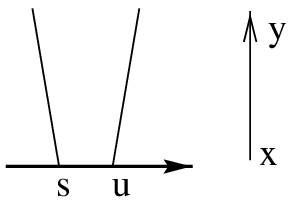}$$
We see that these two labellings only differ by the transposition $(x,s)$, 
thanks to the way the standard orientation of $\sd$ was defined (Subsection
\ref{ori}).
This proves that these two diagrams appear in the relation with different
signs (or on different sides of the equality). \endproof

\section{Anomaly}\label{secano}
\setcounter{equation}{0}

\subsection{The anomaly as a degree of a map}

Recall that the non-degenerated types (a) and (a') faces are defined by
a connected component $A $ of $\g$, whose configurations map all
univalent vertices of $A$ to a line $\D$ directed 
by a tangent vector $s$ to $L$. To study them, we 
need to introduce a few definitions, where we replace $M$ 
by the non-compact line $\R$.
This does not change the fundamental notions of diagrams and configuration
spaces ($\dk(\R)$ is identical to $\dk(J)$ for instance).

\begin{nota} Let $\dw$
be the set of connected elements of $ \dk(\R)$.

For all $s\in \sd $, let $\D$ be the standard map $u\mapsto us$ from $\R$
to the oriented axis of $\rt$ with direction $s$. By an abuse of notation,
$\D$ will also mean this oriented axis. 

For any $\gu\in\dw
$, let $\wx(\gu)$ denote the quotient of $C_\gu(\D)$ by the 
translations and dilations of $\rt$ which preserve the axis $\D$.

Then, let $W(\gu)$ denote the union of the $\wx(\gu)$ for $s$ running
over $\sd $.
\end{nota}

We are going to give these spaces an orientation. First, let us orient
the two-dimensional Lie algebra of the group of translations 
and dilations which preserve $\D$, by giving it a basis: take 
a first vector oriented to the positive translations of the oriented axis
$\D$, and a second vector oriented to the dilations of ratio greater than one.

Now, since the space $C_\gu(\D)$ is oriented by the method of 
Subsection \ref{ori}, 
this orients the quotient space $\wx(\gu)$, and therefore $W(\gu)$
using the standard orientation of the base space $\sd$. 
There is no ambiguity, since all
considered spaces are even-dimensional.

Then we define a map  $\Psi_\gu$ from $W(\gu)$ to 
$(\sd)^{E^v}$, and a map $\Psi'_\gu$ from 
$W'(\gu)=W(\gu)\times (\sd)^{E^a}$ to $(\sd)^{E}$, exactly as before.

Now we are going to see that these spaces can be glued together according to
the same formulas, except that there are no anomalous faces anymore.

\begin{prop}\label{w} For the map $\poids $ defined by the same
formula (\ref{coef}) on $\dw$, the linear combination of spaces
$$W_{n,k}=\sum_{\gu\in\dw
}\poids (\gu) W'(\gu)$$ is completely glued, hence the union $\Psi$ of the
$\Psi'_\gu$ admits a degree on it.
\end{prop}

{\bf Proof of Proposition \ref{w}}

The method is to use the same proofs as those of Sections 4 and 5 in 
this new situation, where there are no equivalents of types (a), (a') and (e)
faces (for the scale is relative and the diagrams are connected): 
we just have to check that the map $\Psi$ defined on the spaces
$W(\gu)$ factorises in the same way
on the faces.
It is very easy to do for the faces $F (\gu,A)$ such that $A\inc T$:
$$F (\gu,A)=W(\gu/A)\times \co(A).$$
The other faces can also be seen as fibered spaces with bases
$W(\gu/A)$ and fiber $W_s(A)$: over each $y\in W(\gu/A)$ we have the fiber
$W_s(A)$ where $y\in W_s(\gu/A)$.

The only new thing to prove is that if an STU' relation involves a 
non-connected diagram, then the faces are degenerated. This can easily be done
by projecting the fiber over $\sd$ to the product of the corresponding 
fibers for the connected components: this product has one dimension less,
corresponding to the relative dilations between the
configurations of the two subdiagrams.
 \endproof
\begin{defin} Let $\ew$ be the set of connected elements of $\ed(\R)$. 
In other words, $\ew$ is the set of isomorphism classes of diagrams in 
the $D^W_{n,2}$ where $n$ runs over $\N$, 
together with a choice of arbitrary orientations.
For all $\gu\in\ew$ we define the 
integral
$$f_\gu=\int_{W(\gu)}\Psi_\gu^*\O.$$
Then we define the anomaly to be 
$$\ano=\sum_{\gu\in\ew}{f_\gu\over2|\gu|}
[\gu]\in \A(\R).$$
\end{defin}

Note that the degree $n$ part $\ano_n $ is half the degree of the map $\Psi$ of the 
above proposition, when $\ank=\A_n(M)$.
Moreover, $\ano_1={[\th]\over 2}$. 
We are going to see that it is also the degree of a map defined on the quotient
of $W_{n,k}$ by a symmetry, with still the same 
coefficients.

\subsection{Symmetry properties of the anomaly}
\def\bg{\bar\gu}

For any labelled oriented diagram $\gu$ with support $\R$, let $\bg$
be the diagram obtained by reversing the injection $i:U\lra\R$ (by the
composition with the involution $u\mapsto -u$ of $\R$),
and transporting the orientations of univalent vertices with $U$ by this
involution:
$$\vep{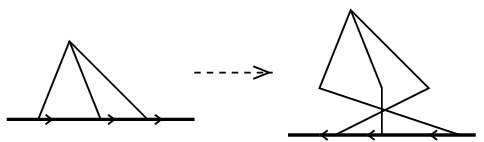}$$
So the order of the vertices is changed.
This defines a linear involution of each $\A_n(\R)$
and an automorphism of $\A(\R)$.

We say that a element of $\A(\R)$ is {\em $S^1$-even} if it is preserved 
by this involution. It is unknown whether all elements of $\A(\R)$ are
$S^1$-even.
Now, let $-\gu$ denote the labelled diagram obtained from $\gu$ by reversing all
its edges and preserving the vertex orientations. We obviously have $[-\gu]=[\gu]$ 
because these diagrams are isomorphic.

We have the two properties:

\begin{lemma}\label{sym1} (i) For all $s\in\sd$, $\gu\in D^W_{n,k}$, there is a negative 
$\Psi$-isomorphism between $W_s(\gu)$ and $W_{-s}(\bar\gu)$. 

(ii)
the anomaly is $S^1$-even ($\bar \ano=\ano$).
\end{lemma}

\begin{lemma}\label{sym2} (i) For all $s\in\sd$ and $\gu\in D^W_{n,k}$, 
there is a 
$\Psi$-isomorphism between $W_s(\gu)$ and $W_{-s}(-\gu)$ with sign $(-1)^n$.

(ii) The anomaly cancels at even degrees.\end{lemma}

{\noindent \bf Proof of Lemma \ref{sym1}}

The $\Psi$-isomorphism that we take from $W_s(\gu)$ to $W_{-s}(\bar\gu)$ is
the identity on configurations as maps from $V$ to $\rt$. The sign comes
from the fact that we quotient the configurations by the action of the
oriented group of dilations and translations along the axis $\D$ which is
reversed. Now, this map preserves the orientation between $W(\gu)$ and
$W(\bar\gu)$ because the antipodal map on the base space $\sd$ also reverses
its orientation. Thus, $f_\gu=f_{\bar \gu}$. \endproof

{\noindent \bf Proof of Lemma \ref{sym2}}

The $\Psi$-isomorphism we take between $W_s(\gu)$ and $W_{-s}(-\gu)$ is
the central symmetry on configurations. Note that this symmetry transports
the orientation of $\D$ and the orientations of univalent vertices.
Thus, the orientation changes come from the geometrical symmetry 
for each trivalent vertex, and the reversal of each edge in the diagram.
So, the sign is $(-1)^n$ because of the formula $\#E-\#T=\deg\gu=n$. Finally,
the antipodal map on the base space and the isomorphism $-\gu\sim\gu$ give 
$f_\gu=(-1)^{n+1}f_{-\gu}=(-1)^{n+1}f_\gu$: $f_\gu$ cancels for even $n$. 
 \endproof\def\dwt{\widetilde D_{n,k}^W}
\def\ddwt{\widetilde D_N^W}

\begin{nota} For $n$ odd, the $\Psi$-isomorphism between 
$W_s(\gu)$ and $W_{-s}(-\gu)$ defines an orientation-preserving involution of $W_{n,k}$. 
We let $\wt_{n,k}$ be the quotient space of $W_{n,k}$ 
by this involution, and let
$\dwt$ be the quotient of $\dw$ by the involution $\gu\mapsto-\gu$.
It lets no diagram fixed because it changes the half-edge of the first
univalent vertex on $\R$. Thus, $$\wt_{n,k}=\sum_{\gu\in\dwt}
\poids (\gu) W'(\gu).$$
\end{nota}

We can note that the image of $\ano_n$ in $\A_n^k(\R)$ is equal to the degree of the map $\Psi$ defined on $\wt_{n,k}$.

\subsection{The labelled generalised diagrams and their configuration spaces}
\label{gen}

Now we are going to define the additional configuration spaces to complete $C'$ in order
to glue all the faces and define the degree of $\Psi$ in the case of a link.

First, identify $M$ with $X\times S^1$ where $X$ is a finite set. 
Then choose a family of smooth maps $\ph_m:\dd\longrightarrow
\sd $ for $m\in X$, where $\dd$ is the disc and each $\ph_m$ 
coincides with the tangent map 
to the link $L$ on the boundary $\{m\}\times S^1$ of $\dd$.
Here we equip $\dd$ with the opposite orientation to 
the usual one. So, the basis (tangent oriented by $S^1$, outward
normal) is direct.

\begin{note} For all $\gu\in\dw$ and $m\in X$, let 
$$\wg= \ph_m^*W(\gu)$$ be the fibered space over $\dd$ which is the union of the
$W_{\ph_m(x)}(\gu)$ for $x\in \dd$, and equip it with the product 
orientation.

In the same way as in Subsection \ref{labd}, we let
$E=\{e_1,\cdots, e_N\}$ be the set of labelled ``edges''
$e_i=\{2i-1,2i\}$. We recall the formula $k=3n-N$. 

Now fix $N$ and let $n$ and $k$ vary for the sets of diagrams $\dw$ and $\dwt$.

Let $\ddw$ denote the union of the $D^W_{n,3n-N}$ for all $n$. 
And let $\ddwt$ denote the union of the $\widetilde D^W_{n,3n-N}$ for all $n$ 
such that $n$ is odd. 
\end{note}

In fact, these unions concern the $n$ such that $2n-1\leq N$, for 
connected diagrams of higher degrees would have more edges.

\begin{defin}\label{gendiag} 
A {\em labelled generalised diagram} with support $M$ is the data of 
a decomposition of $E$ into a disjoint union $$E^v\cup E^a\cup E^W$$
together with:

\b\qua a labelled diagram $\g_0$ with support $M$ such that its set of visible edges is $E^v$,
in the sense of Definition \ref{diag};

\b\qua a set $Y$ of diagrams in $\ddw\sqcup\ddwt$ such that their sets 
of visible edges form a partition of $E^W$, together with a map 
$\gu\mapsto m(\gu)$ from $Y$ to $X$. We let $Y_1=Y\cap\ddw$ and $Y_2=Y\cap
\ddwt$.

The {\em degree} of a labelled generalised diagram is the sum of the degrees of $\g_0$
and all the diagrams $\gu\in Y$.
Let $D^g_{n,3n-N}$ denote the set of these generalised diagrams of degree $n$.
An arbitrary orientation is chosen on each of them.
\end{defin}

Note that $\dl\inc\dg$. Indeed, $\dl$ is the set of generalised diagrams 
$\g\in\dg$ such that $Y=\0$, or equivalently $E^W=\0$.

\begin{defin} Let the {\em configuration space} of a generalised diagram be
$$\c=C_{\g_0}\times\prod_{\gu\in Y_1}\wgg\times 
\prod_{\gu\in Y_2}W(\gu).$$
with the product orientation. This is again completed with the help of the
absent edges:$$\cc=\c\times (\sd)^{E^a}.$$
\end{defin}

The product orientation is well-defined because all spaces are 
even-dimensional. It is reversed when we change the orientation
of a single vertex. The canonical map $\Psi$ extends to the configuration
spaces of generalised diagrams in the obvious way.

\subsection{Gluing the anomalous faces}
\begin{prop}\label{tot}
Let $p$ be a fixed map from $X$ to $\Z$.
We extend $\poids$ to generalised diagrams according to the same formula (\ref{coef}) where now $e_\g=\#(E^v\cup E^W)$ and the images $[\g] \in \A(M) $ of 
the generalised diagrams $\g$ are defined by
\begin{equation}\label{clgen}
[\g]=[\g_0]\prod_{\gu\in Y_1}[\gu]^{(m(\gu))}\times
\prod_{\gu\in Y_2}p_{m(\gu)}[\gu]^{(m(\gu))}.\end{equation}
Then, the canonical map $\Psi$ defined on the linear combination
$$\sum_{\g\in\dg}\poids (\g) \cc$$
admits a degree which belongs to the integral lattice generated by the 
elements $\poids(\g)$ in the space $\ank$ where $\g$ runs over the
principal diagrams in $\dn$.
\end{prop}

\proof

{\noindent \bf Ordinary gluings}

The ordinary gluings of $\g_0$ and all the $\gu\in Y$ are an immediate
consequence of Section \ref{recol} and Proposition \ref{w}.
Here the STU' relation always involves the set $E^a$ of absent edges of the
generalised diagram, given in Definition \ref{gendiag}.
\ppar

{\noindent \bf Anomaly gluings}

Now we are going to glue together the following two types of faces.

\b\qua The types (a) and (a') faces for the subdiagram $\g_0$ of a generalised
diagram; 

\b\qua The faces induced by the boundary of $\dd$ on some $\wgg$.

These faces are glued in pairs in an obvious way, where one face of the first
type is glued to a piece of a face of the second type delimited by the 
positions of the univalent vertices of $\g_0$ on the corresponding component
$\{m\}\times S^1$ of $M$.

Now we have to check that the two orientations are opposite in any case.
The orientation of $\dd$ has been defined by the basis (tangent oriented
by $S^1$, outward normal), whereas for the type
(a) or (a') faces, we must consider the orientation of the Lie algebra
by which we quotiented, when we defined the orientation of the spaces
$W_s(\gu)$: we took the basis (translations in the sense of $s$, growing
dilations). But $s$ represents the generic tangent vector to $L$
with the positive orientation, whereas the growing dilations point inward
to the configuration space.
\ppar

{\noindent \bf Generators} 

By definition, the degree of $\Psi$ is equal to the
count of preimages of a generic point with coefficients $\poids(\g)$. So,
it belongs to the integral lattice generated by the $\poids(\g)$ when $\g$ runs
over the generalised diagrams such that $\Psi(\c)$ has nonempty interior. 
These diagrams are principal, as this necessary condition remains true
generalising Proposition \ref{34} (but the converse is false: a principle
diagram $\gu$ can give an empty interior for a $W(\gu)$).
Then, we may keep the only non-generalised diagrams in this list,
thanks to Formula (\ref{clgen}).
 \endproof

\subsection{Proofs of Propositions \ref{varia} and \ref{mult}}

\begin{nota} For each $m\in X$, define 
the integral $$I_m=\int_{\dd}\ph_m^*\om.$$\end{nota}

\begin{lemma} The integrals on the spaces $\wg$ obey the identity
$$\int_{\wg}\Psi_\gu^*\O= I_m f_\gu.$$
\end{lemma}

This comes from the invariance of the space $W(\g)$ and of the volume form
$\O$ under rotations of $\sd$.
 \endproof

\begin{prop} \label{calcul}
The expression $\deg \Psi$ of Proposition \ref{tot} is 
the projection in $\ank$ of
$$Z \prod_{m\in X}\exp \left((p_m+2I_m) \anom\right).$$
\end{prop}

The proof of this proposition is just a matter of calculation. Let $\egn$ be
the set of isomorphism classes of generalised diagrams in $\dg$. It is
the set of classes of degree $n$ generalised diagrams $\g$ 
such that $e_\g\leq N$.

We have
$$\eqalign{\deg\Psi&=\sum_{\g\in\dg}{(N-e_\g)!\over N!\;2^{e_\g}}
I_L(\g)[\g]\cr
&=\sum_{\g\in\egn}
{I_L(\g_0)[\g_0]\over|\g|}\prod_{\gu\in Y_1} I_{m(\gu)} f_\gu[\gu]^{(m(\gu))}
\times \prod_{\gu\in Y_2}p_{m(\gu)}f_{\gu}[\gu]^{(m(\gu))}.\cr
}$$
\def\bb{\varepsilon_{\gu,m}^\g}
\def\etn{\eta_{\gu,m}^\g}
\def\prodm{\mathop{\prod_{\gu\in\ew}}\limits_{m\in X}}
For any generalised diagram $\g$, $m\in X$ and $\gu\in\ew$, let $\bb$
denote the number of copies of $\gu$ colored by $m$ in the set $Y_1$ of $\g$, 
and let $\etn$ denote its number in the set $Y_2$.

Then, the number of automorphisms of $\g$ is:
$$|\g|=|\g_0|\prodm|\gu|^{\bb}
(\bb)!\;(2|\gu|)^{\etn}(\etn)!$$
So we have
$$\deg\Psi=\!\sum_{\g\in\egn}{I_L(\g_0)[\g_0]\over|\g_0|}
\prodm{(I_m)^{\bb}\over\bb!}\left({f_{\gu} [\gu]^{(m)}
\over|\gu|}\right)^{\bb+\etn}\!\!
{1\over\etn!}\left({p_{m}\over2}\right)^{\etn}.
$$
This is the image in $\A_n^k(M)$ of the degree $n$ part of
$$\left(\sum_{\g_0\in\ed}{I_L(\g_0)[\g_0]\over|\g_0|}\right)
\prod_{m\in X}\exp\left(\sum_{\gu\in\ew}{f_{\gu} [\gu]^{(m)}
\over|\gu|}\left(I_{m}+{p_{m}\over2}\right)\right)$$
which is nothing but
$$Z \prod_{m\in X}\exp \left((p_m+2I_m) \anom\right).\eqno{\qed}$$

\begin{lemma}\label{phiint} For all $m$, $I_L(\th_m)+2I_m\in\Z$.
\end{lemma}

We can restrict ourselves to one fixed component $m$ of $M$, 
thus to the case of a knot.
Apply Proposition \ref{calcul} to the case
$n=N=1$: the degree 1 part of
$Z\exp ((p_m+2I_m) \ano)$ is $I_L(\th){[\th]\over2}+(p_m+2I_m) 
{[\th]\over2}$.
It is the degree of a map where the coefficient is 
$\poids(\th)={[\th]\over2}$.
Thus, $I_L(\th_m)+2I_m$ is an integer. \endproof

\noindent{\bf Proof of Proposition \ref{varia}}

The above constructions work for any map $p$ from $X$ to $\Z$. Now, let
$p_m=-(I_L(\th_m)+2I_m)$. In an isotopy, we let the maps $\ph_m$ 
vary continuously, so the integer $p_m$ is constant for any $m$.

Now let $$Z^0=Z \prod_{m\in X}\exp \left((p_m+2I_m) \anom\right).$$
Its part of each degree is the degree of a map, thus it is constant
in each isotopy class.
We have $p_m+2I_m=-I_L(\th_m)$. We conclude that
$$Z=Z^0 \prod_{m\in X}\exp \left( I_L(\th_m)\anom\right).\eqno{\qed}$$

\noindent{\bf Proof of Proposition \ref{mult}} 

We know that $I_L(\th_m)+2I_m$ is an integer.
If we suppose that $I_L(\th_m)$ is an integer, then we can let $p_m=-2I_m$.
Thus, for any fixed $n$ and $k$, $Z_n^k$ is equal to $\deg\Psi$. The
conclusion comes from Proposition \ref{tot}.
 \endproof

\begin{remark}
It is well-known that for a knot ($I_m=I$, $\th_m=\th$, $\ph_m=\ph$), 
the Gauss integral $I_L(\th)$ 
is some sort of measure of the framing of the knot, that is, the 
linking number of the knot with a parallel curve. Now we have 
proved that $I_L(\th)+2I$ is equal to an integer. So one can wonder 
how to construct geometrically the parallel curve with framing
$I_L(\th)+2I$ from the map $\ph$. Here is the answer:

Consider the fibre bundle over $\dd$, with fiber over $x\in \dd$
the normal plane to the vector $\ph(x)$. Since $\dd$ is homotopically
trivial, there is a unique trivialisation of this bundle up to homotopy.
This gives a homotopy class of trivialisation of the bundle over the boundary
$S^1$ of $\dd$, which is the normal bundle to the knot. Thus, 
it gives our parallel curve up to isotopy.
\end{remark}

\section{Vanishing of the anomaly in degree 3}\label{v3}
\setcounter{equation}{0}

In this section we are going to prove that the degree 3 part of $\ano$ 
vanishes, although
the integrals $f_\gu$ that occur in its expression do not individually cancel.

It is easy to see that the degree 3 diagrams in $\ew$ are $\vep{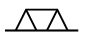}$, $\vep{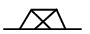}$, and
$\vep{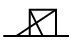}$. They have no automorphisms other than the identity. We
also easily find, when considering their classes in $\A(S^1)$, that
$[\vep{a3}]=-[\vep{a1}]$ and that $[\vep{a2}]=0$. So, what we have to prove
is that $f_{\vep{a1}}=f_{\vep{a3}}$.

Put the following labelling on the diagrams $\vep{a1}$ and $\vep{a3}$:
$$\vep{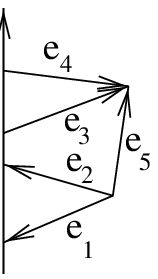}\qquad\vep{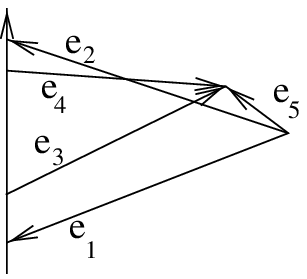}$$
Now we are going to describe the regions $\Psi(C_{\vep{a1}})$ and
$\Psi(C_{\vep{a3}})$ of $(\sd)^E$, $E=\{e_1,\cdots,e_5\}$. They are
the sets of configurations
of $E$ in $\sd$ which verify the two conditions: $$\vep{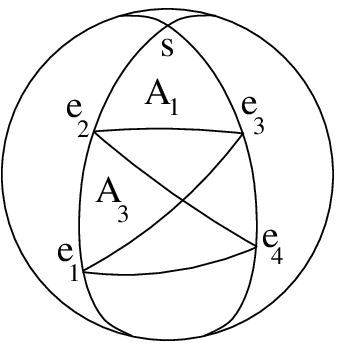}$$
(1)\qua The points $e_1$, $e_2$, $e_3$, $e_4$ form a ``square'',
that is, they obey the four conditions obtained by circular permutations of
this one: the points $e_3$ and $e_4$ are on the same side of the vector 
plane $(e_1, e_2)$ which contains $e_1$ and $e_2$. 

Indeed, considering the vector $s$ orienting $\D$ as the north pole, we can
find it on the picture up to sign as the intersection of $(e_1,e_2)$ and 
$(e_3,e_4)$. We must have $e_1$ and $e_2$ on the same meridian 
(half-circle from $s$), and the same for $e_3$ and $e_4$. This gives two of
the above four conditions.
The two other conditions are equivalent expressions of the fact that, while
the information on the sign of $s$ is given by the orientation of any of these
meridians, respectively from $e_1$ to $e_2$ and from $e_4$ to $e_3$, these
two informations agree.
\smallskip

(2)\qua The point $e_5$ is in the region $A_1$ in the case of $\vep{a1}$, and in
the region $A_3$ in the case of $\vep{a3}$.

Indeed, before fixing $e_5$, the values on the sphere of the $e_i$ for
$i\leq 4$ determine $\D$ and the shape of each of the two halves
of the diagram ($\{e_1,e_2\}$ on the one hand, $\{e_3,e_4\}$ on the other hand)
up to independant translations and dilations preserving $\D$.

Then, in the case of $A_1$ (image of $\vep{a1}$), the vector $e_5$ must be equal to 
$e_2+e_3+\lambda s$ for some $\lambda>0$ and where each $e_i$ is determined only up to a
positive scalar multiplication. Therefore it means that $e_5$ is in the region
$A_1$.

The case of $A_3$ (image of $\vep{a3}$) is defined by the fact that $e_2$ is above $e_4$ and
$e_3$ is above $e_1$. Therefore, $e_5$ is of the form $e_2+e_4-\lambda s$
on the one hand, $e_1+e_3+\lambda' s$ on the other hand. This gives the
region $A_3$.

Note that there are two possible orientations for the square 
$(e_1,e_2,e_3,e_4)$, corresponding
to the two sides of the vector plane $(e_1,e_2)$ where the points $e_3$ and $e_4$ can be.
This splits the configuration spaces into two subsets.
We can restrict the study to one of them, since the other one is obtained from it
by the central
symmetry: we know that this symmetry preserves the integral according to Lemma \ref{sym2}.
On the interior of each of these subsets, $\Psi$ is a diffeomorphism. 

To show the equality of the integrals on these two areas, we shall first
show the equality of signs, then the equality of absolute values.

To prove the equality of signs, let us see the areas $A_1$ and $A_3$ as two regions
of the biangle $B(e_1,e_2,e_3,e_4)$ with one side containing $e_1$ and $e_2$,
and the other side containing $e_3$ and $e_4$. Let us extend the map 
$\Psi$ to a larger ``configuration space'' which contains $C_{\vep{a1}}$ and
$C_{\vep{a3}}$ as open subsets, and in which one can go continuously from one
to the other. 

This larger configuration space is defined as follows. 
It is the quotient by dilations and translations of the set 
of maps from $V$ to $\rt$ such that:

\b\qua  It is injective on each edge.

\b\qua It is not contained in a plane.

\b\qua All univalent vertices are on a straight line.

\b\qua On this line, the respective positions $z_i$ of the ends of 
$e_i$ ($1\leq i\leq 4$) verify $(z_2-z_1)(z_4-z_3)>0$.

Its image by $\Psi$ is the region of $(\sd)^E$ defined by the conditions 
that $e_1$, $e_2$, $e_3$, $e_4$ form a square, and that $e_5$ is inside
the biangle $B(e_1,e_2,e_3,e_4)$.
It is an imbedding, thus its jacobian keeps its sign on each connected 
component.

So, as each connected component meets both $C_{\vep{a1}}$ and
$C_{\vep{a3}}$, we have the same signs on them.

Finally, to prove that the integrals on $C_{\vep{a1}}$ and
$C_{\vep{a3}}$ have the same absolute value, we are going to show that the
system of inequalities delimiting their image (positions of the $e_i$
such that $e_5$ is in $A_1$, resp. $A_3$ relatively to the others) are similar.
Their similarity comes from this remark:

Denoting for three vectors $s_1,s_2,s_3$ of $\rt$, 
$[s_1,s_2,s_3]=\mathop{\rm sign}(\det(s_1,s_2,s_3))$, the condition that 
$(e_1,e_2,e_3,e_4)$ is a square is expressed
by$$[e_1,e_2,e_3]=[e_1,e_2,e_4]=[e_1,e_3,e_4]=[e_2,e_3,e_4].$$
Therefore it is equivalent to saying that $(e_2,e_3,-e_1,-e_4)$ is a square.

Now, the system of inequalities for $A_3$ is that 
$(e_1,e_2,e_3,e_4)$ is a square and $e_5$ is inside this square, in the
triangle limited by the diagonals and one fixed side $(e_1,e_2)$.

The one for $A_1$ is that 
$(e_2,e_3,-e_1,-e_4)$ is a square and $e_5$ is inside this 
square, in the triangle limited by the diagonals and the side $(e_2,e_3)$.
 \endproof

\section{$Z^0$ is a universal Vassiliev invariant}\label{secvassil}
\setcounter{equation}{0}
\def\I{{\cal I}}

We are going to see that the invariant $Z^0$ is
a universal Vassiliev invariant. But there is a little problem:
we have to choose between three possible conventions. Either we study $Z$ as
a Vassiliev invariant of framed links with framings $I_L(\th_m)$ (that can
be supposed to be integers), or we look at $Z^0$ as an invariant of 
zero-framed links (but it is rather artificial as this framing must 
be corrected at each 
desingularisation), or we just consider the common image of all this under
the quotient of $\A(M)$ by the ideal generated by $\th$, as an invariant of
unframed links. We shall consider here the question in terms of the first 
approach.

Let us recall the defining property of universal Vassiliev invariants. 
A map $\I$ from the set of isotopy classes of framed links 
$L:M\hookrightarrow\rt$ to $\A(M)$ is a
universal Vassiliev invariant if: for any
integer $n$ and any
singular oriented link with $n$ double points, that is, an 
immersion $L$ of $M$ into $\rt$ which is an imbedding except for a set $X$
of $n$ points
which have two preimages, the alternate sum of desingularised links
$$\sum_{\zeta:X\lra\{-1,1\}}\left(\prod_i\zeta(i)\right)
\I(L_\zeta)$$
is zero in all degrees $<n$, and is equal in degree $n$ to the 
chord diagram defined by the positions of singular points in $M$.
Here, $L_\zeta$ is the link obtained from $L$ by replacing 
the singularities by the
crossings with signs defined by $\zeta$.

Let us prove this property of $Z$ where the framings are given by the 
$I_L(\th_m)$, considering that the movement of desingularisation is 
infinitesimal.

Consider $Z$ in terms of configured diagrams as described in Subsection \ref{yang}.
The spaces $C$ of the different desingularised links are in correspondence
except at places where some part of graph connects the two strands
together inside the neighbourhood of a desingularisation.
Then, the alternate sum cancels all contributions of the parts of the $C$
in which there is not the presence of such a part of graph near {\em every}
singular point.

The configured diagrams which do contain such a part of 
graph near every
singular point must have degree $\geq n$: indeed, they must have at least one univalent
vertex on each strand of the singular points, so globally at least $2n$ univalent vertices.
The only degree $n$ diagram which has the above property is precisely the
expected chord diagram.

Now, the alternate sum of integrals on these configured chord diagrams
is easily computed: at every neighbourhood of a singular point, the direction
of the chord runs through a hemisphere limited by the tangent plane at
the singularity for each desingularised diagram. The difference gives the
surface of the whole sphere which is equal to $1$ 
according to our definitions.
 \endproof

\appendix\section{Appendix: Changing the map $\poids$}
\setcounter{equation}{0}

In the Appendix, we study the possibility of choosing a map $\poids$ 
different from Formula (\ref{coef}). It is of no usage in the rest of this
article.

\subsection{Weakening the gluing conditions}\label{light}

Finding other maps instead of (\ref{coef}) that still satisfies a gluing to
obtain a result like in Proposition \ref{tot} is something difficult. We shall not see general solutions
but we shall mention some elements of the problem.

The first aspect of the problem is to make a complete list of configuration 
spaces whose image have codimension 0 or 1. 
It seems unlikely that something more could be done than we did, except for
the $W(\gu)$.
As concerns the question of finding possible other configuration spaces
to complete the list, Dylan Thurston \cite{th} proposed to use spaces
constructed in the same way as our $\wg$ where $\dd$ and $\ph$ are
replaced by $C_\th$ and the map $\Psi$ defined on it. But such spaces do
not need to be introduced because they can always be replaced by a 
suitable combination
of the $\wg$ and $W(\gu)$ as we did. I do not know if one
can introduce usefully other spaces.

The second aspect of the problem is to determine which faces are degenerated.
One can sometimes deduce that a face is degenerated from the fact that
it is $\Psi$-isomorphic
to a face of another configuration space with a codimension $\geq2$ image.

We can generalise Lemma \ref{facedil}
to the fact that a type (d) face is degenerated if 
$\exists B\inc A\cap T,\ \#(E'_B\cap E_A)=1$, or if $A$ is not connected.

It should be possible to find more degenerated type (a) faces, for example
the case when at least 3 trivalent vertices are related to 2 univalent 
vertices each. (The corresponding $W(\gu)$ has codimension 1 image).

The third aspect is to make a list of all possible $\Psi$-isomorphisms between 
non-degenerated faces, to enlarge the equivalence relation between faces,
and thus to reduce the number of gluing conditions, so that more maps 
$\poids$ may satisfy them.

We now recall the gluing conditions and give possible ways of weakening them:

(1)\qua The STU' relation (Lemma 5.3): no suggestion.

(2)\qua The IHX' relation (Lemma 5.2), that we can replace by the following 
weaker IHX'' relation: a IHX'' relation is the sum of the 
IHX' relations obtained the following way: we fix the labelled diagram
$\g/e$ where $e\inc T$ is the edge involved, and we let $e$ be replaced by each
absent edge.

(3)\qua The real conditions to glue the type (d) faces are far less constraining 
than to impose $\poids$ to be constant on isomorphism classes of diagrams. 
With one bivalent vertex on $A$ the weaker condition is easy to write as
there is only one way to group in pairs the faces of labelled diagrams for
the $\Psi$-isomorphism we wrote. With more bivalent vertices it is more 
complicated as the same construction gives many possible $\Psi$-isomorphisms;
moreover, all permutations of $E'_A$ may be used. 
In the case of $A\inc T$ with $\#E'_A=3$, the position of
$A$ in the diagram $\g/A$ may be changed.
But these weaker conditions are not easy to make explicit. 

(4)\qua As for the anomaly, we could quotient the configuration spaces of
generalised diagrams by all $\Psi$-isomorphisms used in Lemmas 
\ref{sym1} and \ref{sym2}. Then, a sufficient condition
for the map $\poids$ defined on ordinary
diagrams to be extendable to generalised diagrams in a way which satisfies
the resulting gluing conditions 
is that: for all diagrams obtained by inserting a connected 
$\gu$ into a fixed diagram $\g$ at different places of a fixed 
component of $M$, the expression $$\poids(\g,\gu)+(-1)^{\deg\gu+1}
\poids(\g,-\bar\gu)$$ does not depend on the place where $\gu$ is inserted.
It may not be necessary because some type (a) or (a') faces can be degenerated,
but we cannot easily weaken this as other gluings inside $Y_1$ must be preserved.

We could also wonder if some torsion Vassiliev invariant can be found by this
method, when the map $\poids $ takes values in a torsion $\Z$-module.
But there may also be torsion invariants that cannot be found by this method.
(This is suggested by the fact that the rationality results we can obtain in
degree 2 by this method are not the best ones: see Appendix \ref{deg2}). But
we won't be interested in this problem here.

\subsection{Universality of $Z_n$ in the case of a free module}

Now we are going to see that if the map $\poids$ takes values in a free $\Z$-module 
or in a $\R$-vector space, then the corresponding $\deg\xi$ if it exists, or the 
corresponding integral, comes from $Z_n$.
Thus, the only interest in choosing such other maps $\poids $ is
to obtain better rationality results than those of Theorem 
\ref{mult}.

\begin{prop} If a map $\poids $ from the set of principal elements of 
$\dl$ to some $\R$-vector space
$\cal E$ obeys the STU' relation, then the corresponding expression of
$\deg \Psi$
is of the form $\L(Z_n(K))$ where $\L$ is a linear map from
$\A_n(M)$ to $\ee$ which only depends on $\poids $.
\end{prop}

Let $\enk$ be the set of diagrams 
$\g\in\edm $ which are principal and verify
$u_\g\geq k-1$.

We shall first define a map $\ps$ from $\enk$ to $\ee$. Then,
we shall prove that $\ps$ extends to a linear map $\L$ defined
on $\A_n(M)$, and that this $\L$ is the solution to our problem.

Let us first define $\ps$ from $\enk$ to $\ee$ by:
$$\qs \g,\qua\ps(\g)=|\g|\mathop{\sum_{\g'\in\dl}}\limits_{\g'\sim \g}
\poids (\g').$$
The existence of a linear extension of $\ps$ on $\A_n(M)$ will be deduced
from the following facts:

(1)\qua Any STU' relation implies the corresponding STU relation for $\ps$.

(2)\qua This map $\ps$ is well-defined on chord diagrams and obeys the
four-terms relation.

(3)\qua The STU relation obtained from (1) allows us to express any diagram 
in $\enk$ as a linear combination of chord diagrams.

Indeed, from (2) we deduce that $\ps$ defines a linear map from the space of
degree $n$ chord diagrams with support $M$ modulo the four-terms 
relation. But we know that
this space is identified with $\A_n(M)$ \cite{bn}, 
so $\ps$ defines a linear map $\L$ on $\A_n(M)$.
This map coincides with $\ps$ on $\enk$, thanks  to (1) and (3).

Now, to prove the identity
$$\sum_{\g\in \enk}\poids (\g) I_L(\g)=\sum_{\g\in\edm}{I_L(\g)\over|\g|} \L([\g])$$
we have to check that if a diagram $\g\in\edm $ is not represented in $\enk$,
then $I_L(\g)=0$ or $\L([\g])=0$.
Suppose that $I_L(\g)\not=0$. Then $\g$ is principal, hence
$u_\g<k-1$. Thus, $\g$ can be expressed modulo STU 
as a linear combination of diagrams which have $k-1$ univalent 
vertices and are principal. But since $\L$ vanishes on them,
it also vanishes on $\g$.

Now let us prove the claims (1), (2) and (3).
\ppar

{\bf  Proof of (1)}\qua Consider a STU' relation; let $\g$ be any diagram
in this relation, and let $A$ be the pair of vertices involved.

Note that this STU' relation is defined by the labelled diagram $\g/A$.
We shall consider the sum of this relation
over the set of relabellings of $\g/A$, in other words, the set of the 
$(\Upsilon,\lambda)$, where $\Upsilon$ is a diagram isomorphic to
$\g/A$, and $\lambda$ is an isomorphism from $\g/A$ to $\Upsilon$.

This sum of relations, when applied to one of the three parts of the relation
(the one of $\g$), 
is a sum over the set of relabellings of $\g$, that is, the set of pairs
$(\g',\lambda)$ where $\lambda$ is an isomorphism from $\g$ to $\g'$.
This set is in one-to-one correspondence with the product of the 
set of automorphisms
of $\g$ with the set of diagrams which are isomorphic to $\g$.
So, this sum gives precisely the corresponding term $\ps(\g)$ of the STU relation.

This proves (1). \endproof

{\bf Proofs of (2) and (3)}\qua
Since $k\leq 2n$, $\ps$ is well-defined on the chord diagrams, and also
on the diagrams with only one trivalent vertex. Then the four-terms relation
is a consequence of the STU relation deduced from (1).

The proof of (3) is easy by induction: we just have to note that when 
we use a STU relation to express a diagram in terms of two diagrams 
with one more 
univalent vertex, then these two diagrams are principal 
provided that the first one is. \endproof

This can be generalised to the case when $\poids $ takes values in any free
$\Z$-module $\M$, for $\M$ is canonically imbedded in the real vector
space $\M\otimes\R$.

\subsection{Rationality results in degree 2}\label{r2}\label{deg2}

In this last part of the Appendix, we shall restrict ourselves to the case when $L$ is a knot $K$.

Theorem \ref{mult} gives the same rationality result in degree 2 as the one of 
\cite{le}, that is, a denominator of 48. But the actual denominator is 24,
according to the formula
$$Z_2={I_K(\th)^2\over8}[\th]^2+(v_2-{1\over24})\left[\vep{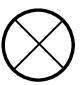}\right]
$$
where $v_2 \in \Z$ is the degree 2 Vassiliev invariant which vanishes
on the unknot $O$ and takes the value 1 on the trefoils. This formula
comes from the Vassiliev universal property of $Z$
(Section \ref{secvassil}), and the computation of 
$$I_O\epfb{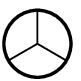}={1\over 8}$$
first made in \cite{gua} (not very difficult).

This actual denominator can be obtained in terms of degree theory 
as Polyak and Viro did \cite{pv}.
Let us recall their remarks, and specify the general context for which
this improvement in degree two is a particular case.

Keep $N=3$, but
consider the map $\poids $ from $D_{2,3}$ to $\Z$ defined by:
$$\poids \epfb{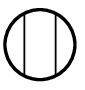}=0, \quad \poids \epfb{cr}={1\over24}$$
independently of the labelling, and $$\poids \epfb{p2}= 0\hbox{ or }-{1\over24}$$
according to the following rule: the labellings of this diagram are
divided into two classes according to the parity of the number
of edges oriented towards the trivalent vertex. Then give the value
$-{1\over24}$ to a class and $0$ to the other.
It is easy to see that the gluing conditions are realised for this
map $\poids $. Therefore, this induces a degree $$\deg \Psi=v_2-{1\over24}.$$
The expression of $Z_2$ suggests that, to find better denominators, it 
would be preferable to manage to cancel the $-{1\over 24}$ term.
This will be done by considering the expression $Z(K)\over Z(O)$.

\begin{prop} The expression $$Z(K)\over Z(O)$$ is equal to the
expression defined the same way as $Z$, with a {\em long knot}:
instead of an imbedding of $S^1$ into $\rt$, we put the knot in the form
of an imbedding
of $\R$ into $\rt$ which coincides with 
the inclusion of a given straight line of $\rt$ outside some compact set.
\end{prop}

We shall not prove this result here as it can be directly deduced from
our next article \cite{p2} on the extension of the configuration space
integral to tangles.

We are only going to show that the result of the integral on this new
shape of knots is an invariant, while the cyclic order on the
univalent vertices of the labelled diagrams is replaced with a total order.

This is expressed in the following lemma:

\begin{lemma} The same conditions as those of Subsection \ref{light} but in
the case of a long knot, imply that $\Psi$ admits a degree.
\end{lemma}

\proof
Let us first define the compactification of the configuration spaces
we use here. We imbed each $\c$ into the compact space 
$$\h\times[-\infty, \1]^U.$$
This compactification has the same general properties as in the case of
compact knots, except that it has new faces which correspond
to the $(f,f')\in \cb$ such that $f'(U)$ meets $\pm\infty$.

Thanks to Remark \ref{discon}, we can suppose that $\g$ is connected.
Now, the fact that $f'(U)$ meets $\pm\infty$ implies that $f_V$ maps all 
univalent vertices to the 
straight line $\lo$. 
We conclude that this face is degenerated by means of the argument of Subsection
\ref{dime}.

Now, let us glue the anomalous faces. Note that the tangent map
d$K: [-\infty, \1]\lra \sd $ verifies
d$K(-\infty)=$ d$K(+\infty)$, so defines a map from
$S^1$ to $\sd $, which is the boundary of some map
$\ph: \dd\lra \sd $. It still has the property that $I_L(\th)+2I$ is an
integer (Lemma \ref{phiint}). (The only difference is that this
integer is now even, whereas it was odd in the case of an ordinary knot.)
 \endproof

With this new problem, if we choose a map $\poids $ which only depends
on the diagram obtained by closing the total order of $U$ into a cyclic 
order, then this gives exactly the same rationality result as the one
we would obtain with this map from the diagrams with a cyclic
order, for the compact knot. So, we can obtain a better rationality result 
for $Z(K)\over Z(O)$ than for $Z(K)$, only by looking for maps 
$\poids $ which do depend on the cut of the cyclic order.

In particular for $n=2$ and $N=3$, take the map $\poids $
which only depends on the orientations of the edges, such that
$$\poids \epfb{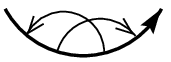}={1\over 6},\qquad \poids \epfb{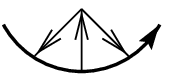}=-{1\over6}$$
and which vanishes on any other diagram. It naturally extends to generalised
diagrams and satisfies all gluings.
The corresponding expression $\deg \Psi$ is equal to $v_2$.

(Actually, no better result can be obtained with another 
map $\poids $ when $N=3$).

\Addresses\recd

\end{document}